\documentclass[11pt]{sat}

\usepackage{amstext}
\usepackage{amsthm}
\usepackage{amsopn}
\usepackage{amsfonts}
\usepackage{amsmath}
\usepackage{amssymb}
\usepackage{amsfonts}

\numberwithin{equation}{section}

\title{Polynomial Approximation and $\omega^r_\varphi  (f,t)$\\
Twenty Years Later}
\def\shorttitle{Polynomial Approximation}

\author{Z. Ditzian}
\def\shortauthor{Z. Ditzian}

\def\versiondate{16 September 2007}

\def\abstracttext{
About twenty years ago the measure of smoothness $\omega  ^r_\varphi
 (f,t)$ was introduced and related to the rate of polynomial
approximation. In this article we survey developments about this and
related concepts since that time.}

\def\MSCnumbers{41A10, 41A17, 41A25, 41A27, 41A30, 41A36, 41A40,
41A50, 41A63, 26A15, 26B35, 26B05, 42C05, 26A51, 26A33, 46E35} 

\def\keywords{Moduli of smoothness, $K$-functionals, realization
functionals, polynomial approximation,
direct and converse inequalities, Bernstein, Jackson, Marchaud,
Nikol'skii and Ul'yanov type inequalities.} 

\def\q{\quad}
\def\qed{\hfill $\square$}
\def\Qbar{\text{\sl Q\kern-.45em{\vrule height.63em width.05em
depth-.033em}}~}
\def\qbar{{{\scriptstyle Q}\kern-.45em{\vrule height.41em width.035em
depth-.03em}}~}
\def\Cbar{\text{\sl C\kern-.35em{\vrule height.63em width.05em
depth-.033em}}~}
\def\cbar{{{\scriptstyle C}\kern-.41em{\vrule height.42em width.035em
depth-.03em}}~}
\def\ibid{\hbox to .5truein{\hrulefill}}
\def\IH{\text{{\rm I}\kern-.13em{\rm H}}}

\def\IR{\text{{\rm I}\kern-.13em{\rm R}}}
\def\pd#1#2{\frac{\partial#1}{\partial#2}}
\def\twoheaddown{\downarrow\kern-0.78em\raise0.25em\hbox{$\downarrow$}}
\def\headtaildown{\downarrow\kern-0.79em\raise 0.5em\hbox{$\ssize\curlyvee$}}

\def\vs{\vskip.3cm}

\def\noi{\noindent}

\def\pr{\prime}

\def\la{\langle}
\def\ra{\rangle}
\def\wt{\widetilde}
\def\wh{\widehat}

\def\os{\overset}
\def\us{\underset}

\def\SD{{\cal D}}

\def\spa{\text{\rm span}}

\def\elra{\hbox to 2in{\rightarrowfill}}
\def\ella{\hbox to 2in{\leftarrowfill}}
\def\hrf{\hbox to 2in{\hrulefill}}
\def\hdotfill{\leaders\hbox to 1em{\hss .\hss}\hfill}
\def\bal{\pmb\alpha  }
\def\bbe{\pmb \beta  }
\def\bh{\pmb h}
\def\bk{\pmb k}

\def\bu{\pmb u}
\def\bv{\pmb v}
\def\bx{\pmb x}
\def\bxi{\pmb \xi  }
\def\by{\pmb y}

\def\floor#1{\lfloor#1\rfloor}

\def\startpagenumber{106}
\def\volumenumber{3}
\def\year{2007}

\newcommand{\beginddoc}{
\begin{document}
\maketitle
\begin{abstract}
\abstracttext
\vskip1pt MSC: \MSCnumbers
\ifx\keywords\empty\else\vskip1pt keywords: \keywords\fi
\end{abstract}
\insert\footins{\scriptsize
\medskip
\baselineskip 8pt
\leftline{Surveys in Approximation Theory}
\leftline{Volume \volumenumber, \year. pp.~\thepage--\pageref{endpage}.}
\leftline{\copyright\ \year\ Surveys in Approximation Theory.}
\leftline{ISSN 1555-578X}
\leftline{All rights of reproduction in any form reserved.}
\smallskip
\par\allowbreak}
\tableofcontents}
\renewcommand\rightmark{\ifodd\thepage{\it \hfill\shorttitle\hfill}\else {\it \hfill\shortauthor\hfill}\fi}
\markboth{{\it \shortauthor}}{{\it \shorttitle}}
\markright{{\it \shorttitle}}
\def\endddoc{\label{endpage}\end{document}}
\date{{\small \versiondate}}
\setlength\oddsidemargin{0pc}
\setlength\evensidemargin{0pc}
\setlength\topmargin{0in}
\setlength\textwidth{6.5in}
\setlength\textheight{8.6in}
\beginddoc

\section{Introduction}\label{Sec1}

It was observed long ago (see \cite{Ni}) that for investigating the rate of
algebraic polynomial approximation the ordinary moduli of smoothness
are not completely satisfactory. For $C[-1,1]$ it was shown that
near the boundary the rate of pointwise approximation was better for
a given degree of smoothness than at other points such as those
further away from the boundary.  The model of the relation between
the ordinary moduli of smoothness and the rate of best trigonometric
approximation (i.e. direct and weak converse inequalities) could not
be followed.  Characterization of the class of functions for which
the rate of best polynomial approximation is prescribed cannot be
described by the ordinary moduli of smoothness.

About twenty years ago the moduli $\omega  ^r_\varphi  (f,t)$ were
introduced (see \cite{Di-To87}) to deal with this problem. There were other
attempts made, the most notable being the works of K.~Ivanov (see
\cite{Iv} for additional references) on the average moduli of smoothness.
 The measure of smoothness $\omega  ^r_\varphi  (f,t)_p$ on $[-1,1]$
(for example) is given by
\begin{equation}\label{Eq1.1}
\omega  ^r_\varphi  (f,t)_p = \;\sup_{\vert  h\vert  \le t}\, \Vert  \Delta  ^r_{h\varphi
}f\Vert  _{L_p[-1,1]}
\end{equation}
where
\begin{equation}\label{Eq1.2}
\Delta  ^r_{h\varphi  }f(x) =\begin{cases}
\overset r{\us{k=0}\sum}\; (-1)^k \binom {r}{k} f\big(x+(\frac r2
-k)h\varphi  (x)\big), \\
{}\q\q\q\q\q \text{\rm if}\q [x-\frac r2\, h\varphi  (x),x+\frac r2\,
h\varphi  (x)]\subset [-1,1]\\
\\
0 \;\q\q\q\q\text{\rm otherwise,}
\end{cases}
\end{equation}
$\varphi  (x)^2 = 1-x^2$ and
$$
\Vert  g\Vert _{L _p[a,b]} = \Big\{\int^b_a \vert  g(x)\vert  ^p
dx\Big\}^{1/p}, \q p<\infty;  \q\q\q \Vert  g\Vert  _{L_\infty[a,b]}  =\us{a\le
x\le b}{\text{\rm ess sup}} \;\vert  g(x)\vert  .
$$

Many properties of $\omega  ^r_\varphi  (f,t)_p$ and
related measures were studied in \cite{Di-To87} as well as the basic relation with
polynomial approximation.  In the last two decades numerous articles
were written using $\omega  ^r_\varphi  (f,t)$ or competing with it.
 In this paper I will give a
survey of what I believe to be the main advances made in the last
twenty years connecting the rate of approximation of functions by
algebraic polynomials with measures of smoothness of these functions.
In \cite{Di-To87} the ``step weight'' function $\varphi  $
was just a function satisfying very mild conditions. Here $\varphi
$ will be a function that is directly used in applications to
approximation and in particular to polynomial approximation and to
some common linear processes.  Unless otherwise specified, when we
write $\omega  ^r_\varphi  (f,t)_p,$ we assume the definition in
(\ref{Eq1.1}) and (\ref{Eq1.2}) on $[-1,1]$ but we will deal also with related
concepts as well as other domains and ``step weights'' $\varphi  .$

We will be discussing relations among different concepts of
smoothness which include $\omega  ^r_\varphi  (f,t),$ various
$K$-functionals, realization functionals, rate of best
approximation, strong converse inequalities as well as the $\tau  $
modulus by Ivanov, moduli given by generalized translations and
others.  Results on the rate of weighted and multivariate polynomial
approximation in relation to various measures of smoothness will also
be described.

The topics are itemized in the Contents (at the
beginning);
however, inevitably some remarks
relating to one topic may appear in a section dedicated to another.
In particular, when a concept or result is introduced in some section, its
relation to items in later sections will be presented in those
sections.

\section{Jackson-type estimates}\label{Sec2}

It is well-known that for $L_p(T),$ where $T$ is the ``circle''
$[-\pi,\pi]$ and
$0<p\le \infty  $,
\begin{equation}\label{Eq2.1}
E^*_n(f)_p \equiv E^*_n(f)_{L_p(T)} \le C\omega  ^r(f,1/n)_{L_p(T)}
\end{equation}
where
\begin{equation}\label{Eq2.2}
E^*_n(f)_{L_p(T)} = \,\text{\rm inf}\, (\Vert  f-T_n\Vert
_{L_p(T)}:T_n\in {\pmb{\cal T}}  _n),
\end{equation}
 ${\pmb{\cal T}}  _n \equiv\;\text{\rm span}\, \{e^{ikx}:\vert k\vert  <
n\}$ is the set of trigonometric polynomials of degree less than $n$
for $n=1,2,\dots$,
and
\begin{equation}\label{Eq2.3}
\begin{aligned}
\omega  ^r(f,t)_{L_p(T)} &= \us{\vert  h\vert\le t}\sup\;\Vert
\Delta  ^r_h f\Vert  _{L_p(T)}, \\
\Delta  ^r_h f(x) &=\sum^r_{k=0}
(-1)^k\binom{r}{k} f\Big(x+\big(\frac r2-k\big)h\Big),
\end{aligned}
\end{equation}
for $r=0,1,2,\dots$ are the ordinary $L_p$ moduli of smoothness.
In fact (\ref{Eq2.1}) is valid also if $L_p(T)$ is replaced by a
Banach space $B$ of functions on $T$
satisfying
\begin{equation}\label{Eq2.4}
\Vert  f(\cdot + a)\Vert  _B = \Vert  f(\cdot)\Vert  _B \q\q \forall
\; a \in \IR
\end{equation}
and
\begin{equation}\label{Eq2.5}
\Vert  f(\cdot + h) -f(\cdot)\Vert  _B = o(1), \q h\to 0 ;
\end{equation}
that is,
$$
E^*_n(f)_B = \;\text{\rm inf}\, (\Vert  f-T_n\Vert  _B:T_n\in
{\pmb{\cal T}}
_n)\le C\omega  ^r(f,1/n)_B \eqno(2.1)^\pr
$$
where $E^*_n(f)_B$ and $\omega  ^r(f,1/n)_B$ are given by
(\ref{Eq2.2}) and (\ref{Eq2.3}) with $B$ replacing $L_p(T).$  (See
Appendix for a proof of $(2.1)^\pr.)$

For $L_p[-1,1],$ $1\le p\le \infty , $ it was proved in \cite[Theorem
7.2.1]{Di-To87} that
\begin{equation}\label{Eq2.6}
E_n(f)_p\equiv E_n(f)_{L_p[-1,1]} \le C\omega  ^r_\varphi
(f,1/n)_{L_p[-1,1]}
\end{equation}
where
\begin{equation}\label{Eq2.7}
E_n(f)_{L_p[-1,1]} = \;\text{\rm inf}\, (\Vert  f-P_n\Vert
_{L_p[-1,1]} :P_n\in {\Pi}_n),
\end{equation}
 ${\Pi}_n \equiv\;\text{\rm span}\, (1,x,\dots,x^{n-1})$ is the
set of algebraic polynomials of degree at most $n-1$
and $\omega
^r_\varphi  (f,t)_p$ is given by (\ref{Eq1.1}) and (\ref{Eq1.2}).

DeVore, Leviatan and Yu \cite[Theorem 1.1]{De-Le-Yu} showed that
(\ref{Eq2.6}) is valid for $0<p<1$ as well.  The method of their
proof uses a Whitney-type estimate by polynomials of degree $r-1$
and ``patching'' them up by polynomials of degree $n$ that form a
partition of unity, (see also the remark in \cite[p.~74]{Di-Hr-Iv}
about the necessity of Lemma~5.2 there for their proof). This
type of argument is used in \cite{De-Lo} to prove the result for
$1\le p\le \infty  $ as well.

For $L_p[-1,1]$ and other spaces a Jackson-type estimate using a
measure of smoothness given by a $K$-functional which is not always
equivalent to $\omega  ^r_\varphi  (f,t)$ but is still optimal (in
the same sense) will be discussed in Section \ref{Sec4}.  However,
(\ref{Eq2.6}) was not extended to a form which follows $(2.1)^\pr.$
That is, we do not have (\ref{Eq2.6}) with $B$ (satisfying some
general conditions) replacing $L_p[-1,1].$

It was proved by M. Timan [Ti,M,58] that for trigonometric
polynomials a sharper (than (\ref{Eq2.1})) Jackson-type inequality holds,
i.e. for $E^*_k(f)_p$ of (\ref{Eq2.2})
\begin{equation}\label{Eq2.8}
\begin{gathered}
n^{-r}\Big\{\sum^n_{k=0} k^{sr-1} E^*_k(f)^s_p\Big\}^{1/s} \le
C(r,s,p)\omega  ^r(f,n^{-1})_p, \\
s=\max(p,2), \q 1<p<\infty  .
\end{gathered}
\end{equation}
This result, which is best possible for $1<p<\infty  ,$ has rarely
been cited in literature in the English language
and I could find it only in a text by Trigub and Belinsky
\cite[p.~191, 4.8.8]{Tr-Be} (and there without proof and with
$n^{-r}$ missing on the left of ({\ref{Eq2.8})).

Recently, an analogue of this result was proved in \cite{Da-Di-Ti},
that is
\begin{equation}\label{Eq2.9}
\begin{gathered}
n^{-r}\Big\{\sum^n_{k=r} k^{sr-1}E_k(f)^s_p\Big\}^{1/s}
\le C(r,s,p)\omega  ^r_\varphi  (f,n^{-1})_p,\\
s = \max (p,2), \q 1<p<\infty
\end{gathered}
\end{equation}
where $\omega  ^r_\varphi  (f,t)_p$ and $E_k(f)_p$ are given in
(\ref{Eq1.1}) and (\ref{Eq2.7}).

We note that in \cite{Da-Di-Ti} (\ref{Eq2.9}) is just one of many
related formulae and the treatment in \cite{Da-Di-Ti} uses best
approximation by various systems of functions and various measures
of smoothness.

We also note that for $1<p<\infty  $ (\ref{Eq2.9}) was shown in
\cite{Da-Di-Ti}
to be
equivalent to
\begin{equation}\label{Eq2.10}
t^r\Big\{\int^{1/2}_t\, \frac{\omega  ^{r+1}_\varphi
(f,u)^s_p}{u^{sr+1}}\;du\Big\}^{1/s}\le C\omega  ^r_\varphi  (f,t)_p,
\end{equation}
for $1<p<\infty  $ and $s=\max\,(p,2).$

Examples were given in \cite[Section 10]{Da-Di-Ti} to show that (as
far as $s$ is concerned) the inequalities (\ref{Eq2.9}) and
(\ref{Eq2.10}) are optimal for $1<p<\infty  .$

The inequality
(\ref{Eq2.10}) is sharper than the inequality
\begin{equation}\label{Eq2.11}
\omega  ^{r+1}_\varphi  (f,t)_p\le C\omega  ^r_\varphi  (f,t)_p
\end{equation}
for the range $1<p<\infty  .$  The inequality (\ref{Eq2.11}),
however, is valid for the bigger range $0<p\le \infty$ (see \cite
[Chapter 7]{Di-To87} and \cite{Di-Hr-Iv}).

\section{$K$-functionals}\label{Sec3}
As an alternative to $\omega  ^r_\varphi  (f,t)$ one can measure
smoothness using $K\text{-functionals}.$

It was shown in \cite[Theorem 2.1.1]{Di-To87} (not just for the case
$\varphi  (x)^2=1-x^2)$ that
\begin{equation}\label{Eq3.1}
K_{r,\varphi  }(f,t^r)_p\approx \omega  ^r_\varphi  (f,t)_p, \q 1\le
p\le \infty,
\end{equation}
that is
\begin{equation}\label{Eq3.2}
C^{-1}K_{r,\varphi  }(f,t^r)_p\le \omega  ^r_\varphi  (f,t)_p\le
CK_{r,\varphi  }(f,t^r)_p, \q 1\le p\le\infty  ,
\end{equation}
where
\begin{equation}\label{Eq3.3}
K_{r,\varphi  }(f,t^r)_p = \inf\big(\big\Vert  f-g\big\Vert
_{L_p[-1,1]} + t^r\big\Vert  \varphi  ^rg^{(r)}\big\Vert
_{L_p[-1,1]}: g,\dots, g^{(r-1)}\in A.C._{\ell\text{\rm oc}}\big).
\end{equation}
In fact, it is known that in (\ref{Eq3.3}) $g,\dots, g^{(r-1)}\in
A.C._{\ell\text{\rm oc}}$ can be further restricted using instead $g\in C^r[-1,1]$ or
even $g\in C^\infty  [-1,1]$ without any effect on (\ref{Eq3.1}).
One could have observed that $g\in C^r[-1,1]$ is sufficient
already from the proof in \cite{Di-To87}.  That it is sufficient to
consider $g$ in the class $C^\infty  [-1,1]$ follows from the
realization results mentioned in Section \ref{Sec5}.  We note also
that for $p=\infty  $ the result is of significance only when $f\in
C[-1,1]$ as otherwise neither side of (\ref{Eq3.1}) is small when
$t$ is.

For the well-studied analogue on the circle $T$ one has
\begin{equation}\label{Eq3.4}
\omega  ^r(f,t)_B\approx \inf \big(\big\Vert  f-g\big\Vert  _B +
t^r\big\Vert  g^{(r)}\big\Vert  _B:\,g^{(r)}\in B\big) =
K_r(f,t^r)_B
\end{equation}
where $B$ is any Banach space of functions on $T$ in which
translations are  continuous isometries, that is  translations
satisfy (\ref{Eq2.5}) and (\ref{Eq2.4}) respectively.
The notation $g^{(r)} \in B$ means that the $r\text{\rm -th}$ derivative in
${\cal S}^\pr$ (the space of tempered distributions) is in $B.$

We will often use the notation $A(t)\approx B(t)$ and, following
(\ref{Eq3.2}), we mean $C^{-1}B(t)\le A(t)\le CB(t)$ for all relevant
$t.$

We do not have
\begin{equation}\label{Eq3.5}
K_{r,\varphi  }(f,t^r)_B\approx \omega  ^r_\varphi  (f,t)_B
\end{equation}
where $\Vert  f\Vert  _{L_p[-1,1]}$ is replaced by $\Vert  f\Vert
_B$ for a ``general'' Banach space on $[-1,1].$

For an Orlicz space
of functions on $[-1,1]$ this was done in \cite{Wa} in his thesis in
Chinese (and I believe also earlier).  Not being able to read that
work, I cannot describe it.  I learned about it from its
extension to the multivariate situation in \cite{Zh-Ca-Xu} where the
univariate case is taken for granted.

In the next section different
but related $K\text{-functionals}$ will be described for which the
treatment for various spaces is given.

For $L_p[-1,1]$ when $0<p<1$ it was shown in \cite{Di-Hr-Iv} that for
all $f$ in $L_p[-1,1],$ \; $0<p<1,$
\begin{equation}\label{Eq3.7}
K_{r,\varphi  }(f,t^r)_p = 0
\end{equation}
where $K_{r,\varphi  }(f,t^r)_p$ is defined by (\ref{Eq3.3}) with
the quasinorm $\Vert  \cdot\Vert  _{L_p[-1,1]}.$ The proof in
\cite{Di-Hr-Iv} is univariate and local and applies to the circle
$T$ as well, that is
$$
f\in L_p
(T) \q\text{\rm implies}\q
K_r(f,t^r)_p = 0 \q\text{\rm for}\q 0<p<1. \eqno(3.6)^\pr
$$

The identity (\ref{Eq3.7}) implies that we cannot have (\ref{Eq3.2})
for $0<p<1$ as $\omega  ^r_\varphi  (f,t)_p$ is not always zero.
(Clearly, $\vert  x\vert  \in L_p[-1,1]$ and $\omega  _\varphi
(f,t)_p\equiv \omega  ^1_\varphi  (f,t)_p\ne 0.)$ Even before
(\ref{Eq3.7}) was proved, it was
clear that $\omega  ^r_\varphi  (f,t)_p$ cannot be equivalent to
$K_{r,\varphi  }(f,t^r)_p$ when $0<p<1,$ as the saturation rate of
$\omega  ^r_\varphi  (f,t)_p$ is $O(t^{r-1+\frac 1p}
)$ for that
range and $K_{r,\varphi  }(f,t^r)_p$ as a $K$-functional
cannot tend to zero at a rate
faster than $t^r$ unless it equals $0.$

\section{$K$-functionals (second approach)}\label{Sec4}

For a Banach space $B$ of functions on domain $\SD$ and a
differential operator $P_r(D)$ of degree $r$ we define the
$K\text{\rm -functional}$
\begin{equation}\label{Eq4.1}
K_{rm}\big(f,P_r(D)^m,t^{rm}\big)_B\equiv \inf \big(\Vert  f-g\Vert
_B + t^{rm}\Vert  P_r(D)^mg\Vert  _B: P_r(D)^mg\in B\big).
\end{equation}
One can assume $P_r(D)^mg $ is defined as a distributional derivative,
and in most cases we deal with we may assume $g\in C^{rm}(\SD)$
without changing the asymptotic behaviour of
$K_{rm}\big(f,P_r(D)^m,t^{km}\big)_B$ given in (\ref{Eq4.1}). The
$K\text{\rm -functional}$ $K_{r,\varphi  }(f,t^r)_p$ of (\ref{Eq3.3}) is
$K_r\big(f,P_r(D),t^r\big)_p$ with $P_r(D) = \varphi
^r\big(\frac{d}{dx}\big)^r$ on $L_p[-1,1].$  In relation to
polynomials on $[-1,1]$ it is natural to study the $K\text{\rm
-functional}$
given in (\ref{Eq4.1}) with $P_2(D) =\frac{d}{dx}\,
(1-x^2)\,\frac{d}{dx}\,.$ It was essentially shown in \cite[Theorem
5.1]{Ch-Di94}, using a maximal function estimate, that
\begin{equation}\label{Eq4.2}
K_{2,\varphi  }(f,t^{2})_{L_p[-1,1]} \le
CK_{2}\Big(f,\,\frac{d}{dx}\,(1-x^2)\,\frac{d}{dx}, t^2\Big)_{L_p[-1,1]}\q\text{\rm for}\q 1<p\le \infty  .
\end{equation}
It follows from \cite[Chapter 9, 135-6]{Di-To87} which uses
 the Hardy inequality, that for $1\le p <\infty  $
\begin{equation}\label{Eq4.3}
K_{2r}\Big(f,\big(\frac{d}{dx}\,(1-x^2)\,\frac{d}{dx}\,\big)^r, t^{2r}\Big)
_{L_p[-1,1]} \le CK_{2r,\varphi  }(f,t^{2r})_p + t^{2r}
E_1 (f)_p\,.
\end{equation}
It can easily be deduced from \cite[Theorem 7.1]{Da-Di05} that
for $1<p<\infty  $
\begin{equation}\label{Eq4.4}
K_{2r}\Big(f,\big(\frac{d}{dx}\,
(1-x^2)\,\frac{d}{dx}\big)^r,t^{2r}\Big)_{L_p[-1,1]} \approx
K_{2r,\varphi  }(f,t^{2r})_p + t^{2r} E_1(f)_p\,.
\end{equation}
For $p=1$ and $p=\infty  $ (\ref{Eq4.4}) does not hold, as shown in
\cite[Remark 7.9, p.88]{Da-Di05}.

We observe that for $r=1$ (\ref{Eq4.4}) is a corollary of
(\ref{Eq4.2}) and (\ref{Eq4.3}), whose proof is more elementary.
(It does not use the Muckenhoupt transplantation theorem nor the
H\"ormander-type multiplier theorem used for the proof of
\cite[Theorem7.1]{Da-Di05}.)  It would be nice if we had a proof for
(\ref{Eq4.2}) with $2$ replaced by $2r$ and could deduce
(\ref{Eq4.4}) directly from it and (\ref{Eq4.3}).

For an orthonormal sequence of functions $\{\varphi  _n\}$ on some
set $D$ the
Ces\`aro summability of order $\ell$ is given by
\begin{equation}\label{Eq4.5}
C^\ell_n(f,x) =\sum^n_{k=0} \,\Big(1-\frac{k}{n+1}\Big)\cdots
\Big(1-\frac{k}{n+\ell}\Big)P_k(f,x)
\end{equation}
where the $(L_2$ type) projection $P_kf$ is given by
\begin{equation}\label{Eq4.6}
P_k(f,x) = \varphi  _k(x)\int_D \varphi  _k(y)f(y)dy.
\end{equation}
Here $D = [-1,1]$ and $\varphi  _k(x)$ are the eigenvectors of
$\frac{d}{dx}\, (1-x^2)\,\frac{d}{dx}$ satisfying
\begin{equation}\label{Eq4.7}
\frac{d}{dx}\, (1-x^2)\,\frac{d}{dx}\,\varphi  _k (x) =
-k(k+1)\varphi  _k(x), \q \int^1_{-1} \varphi  _k(x)\varphi  _\ell
(x)dx = \begin{cases} 0, &k\ne \ell,\\ 1, &k=\ell.\end{cases}
\end{equation}
In later sections we deal with weights in (\ref{Eq4.6}) and
(\ref{Eq4.7}) when we discuss progress made for measures of
smoothness and polynomial approximation in weighted $L_p$ and in
other related Banach spaces. Furthermore, it will be crucial to
examine (\ref{Eq4.5}) when the projection is on a finite dimensional
orthonormal space which is needed for the multivariate situation
(and has the precedent of projection on $\text{span}\,(\sin
kx,\,\cos kx)).$

The Legendre operator $\frac{d}{dx}\, (1-x^2)\,\frac{d}{dx}$ has as
eigenvectors the Legendre orthogonal polynomials.  It was shown
in \cite[Theorem 4.1 and (6.13)]{Ch-Di97} and \cite{Di98} that for $B$ a Banach space of
functions on $[-1,1]$ for which
\begin{equation}\label{Eq4.8}
\Vert  C^\ell _n (f,\cdot)\Vert  _B\le C\Vert  f\Vert  _B
\end{equation}
is satisfied for some $\ell,$ one has
\begin{equation}\label{Eq4.9}
E_n(f) _B =\underset{P\in \Pi_n}\inf\;\Vert  f-P\Vert  _B \le
CK_{2r}\Big(f,\big(\frac{d}{dx}\,(1-x^2)\,\frac{d}{dx}\big)^r,t^{2r}
\Big)_B.
\end{equation}

It is known that $B= L_p[-1,1]$ satisfies (\ref{Eq4.8}) (see for
discussion and references of more general results \cite[Theorem A,
page 190]{Ch-Di97})
and perhaps this
should be an incentive to investigate for which class of Banach
spaces (\ref{Eq4.8}) is valid (with respect to eigenfunctions of
$\frac{d}{dx}\, (1-x^2)\,\frac{d}{dx}),$ and hence imply
(\ref{Eq4.9}) which is a Jackson-type result for a different measure
of smoothness.

For $\alpha  >0,$ the operator
$\big(-\;\frac{d}{dx}\, (1-x^2)\,\frac{d}{dx}\big)^\alpha  g$ is
defined by
\begin{equation}\label{Eq4.10}
\Big(-\;\frac{d}{dx}\,(1-x^2)\;\frac{d}{dx}\Big)^\alpha  g\sim
\sum^\infty  _{k=1} \lambda  (k)^\alpha  P_k g, \q \lambda  (k) =
k(k+1)
\end{equation}

and we say $\big(-\,\frac{d}{dx}
\,(1-x^2)\,\frac{d}{dx}\big)^\alpha  g\in B$ if there exists a
function $G_\alpha  \in B$ which satisfies $P_k G_\alpha  = \lambda
(k)^\alpha  P_kg.$
We may define the $K\text{\rm -functional}$ (see \cite[p.~324]{Di98}) by
\begin{equation}\label{Eq4.11}
K_{2\alpha  }\Big(f,\big(-\,\frac{d}{dx}\,
(1-x^2)\;\frac{d}{dx  }\big)^\alpha  ,t^{2\alpha  }\Big)_B =
\inf\, \Big(\Vert  f-g\Vert  _B + t^{2\alpha  }\Big\Vert
\big(-\,\frac{d}{dx}\,(1-x^2)\;\frac{d}{dx}\big)^\alpha  g\Big\Vert
_B\Big)
\end{equation}
where the infimum is taken on $g$ such that $g\in B$ and
$\big(-\frac{d}{dx}(1-x^2)\frac{d}{dx} g\big)^\alpha  \in
B.$  For integer $\alpha=r$, (\ref{Eq4.11}) and (\ref{Eq4.1}) with $P_r(D)
=\big(\frac{d}{dx}\, (1-x^2)\;\frac{d}{dx}\big)^r$ are the same
concept.  In \cite[Theorem 6.1]{Di98} it was shown for $B$ such
that (\ref{Eq4.8}) is satisfied that
\begin{equation}\label{Eq4.12}
E_n (f)_B \le CK_{2\alpha
}\Big(f,\big(-\,\frac{d}{dx}\,(1-x^2)\,\frac{d}{dx}\big)^\alpha
,1/n^{2\alpha  }\Big)_B.
\end{equation}

\section{Realization}\label{Sec5}
Realization functionals were introduced by Hristov and Ivanov
\cite{Hr-Iv} in order to characterize $K$-functionals.  As it
happened, this concept gained in usefulness when it was observed
that certain $K$-functionals are always equal to zero for $0<p<1$
(see (\ref{Eq3.7}) or $(3.6)^\pr),$ and one needs an expression that
will replace the $K$-functional and will yield a meaningful measure
of smoothness for all $0<p\le \infty  .$  Realization functionals
were shown in \cite{Di-Hr-Iv} to be such a concept.  It is a
mistake, however, to think that realizations are useful only for
$0<p<1.$  Many articles, starting with \cite{Hr-Iv}, utilized
properties of realizations for various applications.  We will
present here realization-functionals that are measures of smoothness
related to polynomial approximation and $\omega  ^r_\varphi
(f,t)_p\,.$

The most common realization related to $\omega  ^r_\varphi  (f,t)_p$
is
\begin{equation}\label{Eq5.1}
R_{r,\varphi  }(f,n^{-r})_p = \Vert  f-P_n\Vert  _{L_p[-1,1]} +
n^{-r}\Vert  \varphi  ^r P^{(r)}_n\Vert  _{L_p[-1,1]}
\end{equation}
where $P_n\in \Pi_n$ is the best polynomial approximant from
$\Pi_n$ to $f$ in $L_p,$ that is
\begin{equation}\label{Eq5.2}
E_n(f)_p = \us{P\in \Pi_n}{\text{\rm inf}} \; \Vert  f-P\Vert
_{L_p[-1,1]} = \Vert  f-P_n\Vert  _{L_p[-1,1]},\qquad P_n\in \Pi_n
\end{equation}
or a near best polynomial approximant
\begin{equation}\label{Eq5.3}
\Vert  f-P_n\Vert_{L_p[-1,1]} \le AE_n(f)_p,\qquad P_n\in \Pi_n
\end{equation}
with $A$ independent of $n$ and $f.$  Sometimes it is convenient to
use $P_n$ as a polynomial of degree $mn$ which satisfies
(\ref{Eq5.3}).  A particularly convenient polynomial of this nature
for $1\le p\le \infty  $ is the de~la~Vall\'ee Poussin-type operator
on $f$ given by
\begin{equation}\label{Eq5.4}
\eta  _n f =\sum^{2n}_{k=0}\eta  \Big(\frac kn\Big)P_kf
\end{equation}
where $P_kf$ is given by (\ref{Eq4.6}) and (\ref{Eq4.7}), $\eta
(y)\in C^\infty  [0,\infty  ),$ $\eta  (y)=1$ for $y\le 1$ and $\eta
 (y) =0$ for $y\ge 2.$  Clearly, $\eta  _n f\in \Pi_{2n},$ $\eta
_nP=P$ for $P\in \Pi_n,$ and it is known that $\Vert  \eta  _n
f\Vert  _p \le C\Vert  f\Vert  _p$ for $1\le p\le \infty  .$
The inequality $\Vert \eta  _nf\Vert  _B \le C\Vert  f\Vert  _B$ for
$B= L_p[-1,1]$ (and in fact for any $B$ satisfying (\ref{Eq4.8}))
follows the same method used in \cite[p.~192]{Ch-Di97} and
\cite[p.~326--327]{Di98} (using the Abel tranformation and $P_k f =
\os\leftarrow\Delta  ^{\ell+1} \begin{binom}{k+\ell}{\ell}\end{binom} C^\ell_k f$
where $\os\leftarrow\Delta  a_k = a_k - a_{k-1}$ and
$\os\leftarrow\Delta  ^m a_k = \os\leftarrow\Delta
\,(\os\leftarrow\Delta  ^{m-1}a_k)).$  Other de~la~Vall\'ee
Poussin-type operators (or delayed means) were also used for
realizations (see for instance \cite{Ch-Di97} and \cite{Di98}.
The
advantage of using a de~la~Vall\'ee Poussin-type operator (in some
form) over using the best approximant is threefold: it is given by a
linear operator, it is often independent of $1\le p\le
\infty , $ and it commutes with the differential operator
$\frac{d}{dx}\,(1-x^2)\;\frac{d}{dx}\,.$

We can also define (as was originally done)
\begin{equation}\label{Eq5.5}
R^*_{r,\varphi  }(f,n^{-r})_p = \us{P\in \Pi_n}{\text{\rm
inf}}\;\big(\Vert  f-P\Vert  _{L_p[-1,1]} + n^{-r}\Vert  \varphi  ^r
P^{(r)}\Vert  _{L_p[-1,1]}\big).
\end{equation}

It is known and easy to show that (\ref{Eq5.1}) with $P_n$ of
(\ref{Eq5.2}) or (\ref{Eq5.3}) and (\ref{Eq5.5}) are equivalent for
$0<p\le\infty  ,$ that is, $R^*_{r,\varphi  }(f,n^{-1})_p\approx
R_{r,\varphi  }(f,n^{-1})_p.$  If we use $\eta  _nf$ of
(\ref{Eq5.4}) in (\ref{Eq5.1}) for $P_n,$ the equivalence holds only
for $1\le p\le \infty  .$ ((\ref{Eq5.4}) is not defined for $0<p<1.)$

It was proved in \cite{Di-Hr-Iv} that
\begin{equation}\label{Eq5.6}
R^*_{r,\varphi  }(f,n^{-r})_p \approx \omega  ^r_\varphi
(f,n^{-1})_p
\end{equation}
for $0<p\le \infty,  $ and hence $R_{r,\varphi  }(f,n^{-r})_p\approx
\omega  ^r_\varphi  (f,n^{-1})_p$ for $0<p\le \infty  $ if $P_n$ is
given by (\ref{Eq5.2}) or (\ref{Eq5.3}), and for $1\le p\le \infty
$ if for $P_n$ we write $\eta  _n f$ given in (\ref{Eq5.4}).

We note (see \cite{Di-Hr-Iv}) that an analogous result to
(\ref{Eq5.6}) is known for $L_p(T)$ where $T_n,$ an $n\text{\rm
-th}$ degree trigonometric polynomial, replaces $P_n,$ and $\omega
^r(f,t)_p$ replaces $\omega  ^r_\varphi  (f,t)_p.$  The equivalence
(\ref{Eq5.6}) was also extended to other realizations and measures
of smoothness.

For $L_p[-1,1],$ $1\le p\le\infty  $ and other Banach spaces some
sequences of linear operators $A_nf$ other than $\eta  _n f$ given
in (\ref{Eq5.4}) were used for defining the realization
\begin{equation}\label{Eq5.7}
\wt R_{r,\varphi  }(f,n^{-r})_{L_p[-1,1]} = \Vert  f-A_nf\Vert
_{L_p[-1,1]} + n^{-r}\Vert  \varphi  ^r(A_nf)^{(r)}\Vert
_{L_p[-1,1]}
\end{equation}
(see, for instance, \cite{Ch-Di97} and \cite{Di98}).

Of course for $\wt R_{r,\varphi  }(f,n^{-r})_{L_p[-1,1]},$ $A_n$
may depend on $r.$  We will encounter some natural expressions of
the form (\ref{Eq5.7}) in this survey.  In most situations here when
dealing with (\ref{Eq5.7}) either the choice (\ref{Eq5.1}) where
$P_n = \eta  _nf$
with
$\eta  _nf$ of (\ref{Eq5.4})
(which is a near best approximant)
is more useful or
we have a linear
approximation process
$A_nf$ which satisfies a relation with $\omega  ^r_\varphi  (f,t)_p$ that
is superior to $\wt R_{r,\varphi  }(f)_p\approx \omega  ^r_\varphi
(f,n^{-r})_p$ (see Section~8).  The conditions that $P_n$ satisfies,
(\ref{Eq5.2}), (\ref{Eq5.3}) or (\ref{Eq5.4}), are independent of $r$
and this fact has proved useful in many applications. We note that in
the expression $R^*_{r,\varphi  }(f,n^{-r})_p,$ $P_n$ depends on $r$
and hence in applications it is sometimes more advantageous to use
the equivalent form $R_{r,\varphi}  (f,n^{-1})_p.$

For a general Banach space $B$ on $[-1,1]$
it is convenient to deal with
\begin{equation}\label{Eq5.8}
R_{2\alpha  }\big(f,P(D)^\alpha  ,n^{-2\alpha  }\big)_B = \Vert  f-P_n\Vert
 _B +\frac{1}{n^{2\alpha  }}\;\big\Vert  \big(P(D)\big)^\alpha
P_n\big\Vert  _B
\end{equation}
where $P(D)$ is the Legendre operator $P(D) = -\,\frac{d}{dx}\,
(1-x^2)\, \frac{d}{dx}\,,$ $P_n$ is given by (\ref{Eq5.2}),
(\ref{Eq5.3}) or (\ref{Eq5.4}), and $\big(P(D)\big)^\alpha  $ is
given by (\ref{Eq4.9}).
We have (see \cite{Da-Di05})
\begin{equation}\label{Eq5.9}
R_{2r}\big(f,P(D)^r,n^{-2r}\big)_p\approx R_{2r,\varphi
}(f,n^{-2r})_p + n^{-2r}\,E_1(f)_p\q\text{\rm for}\q
1<p<\infty  .
\end{equation}
However, (\ref{Eq5.9}) is not valid for $p=1$ and $p=\infty  $ since
\begin{equation}\label{Eq5.10}
R_{2r}\big(f,P(D)^r,n^{-2r}\big)_p\approx
K_{2r}\big(f,P(D)^r,n^{-2r}\big)_p, \q 1\le p\le \infty  ,
\end{equation}
and also since (\ref{Eq4.4}) is not valid for $p=1$ and $p=\infty  .$ In
fact, for any Banach space $B$ for which (\ref{Eq4.8}) is satisfied
we have
\begin{equation}\label{Eq5.11}
R_{2\alpha  }\big(f,P(D)^\alpha  ,n^{-2\alpha  }\big)_B\approx
K_{2\alpha  }\big(f,P(D)^\alpha  ,n^{-2\alpha  }\big)_B
\end{equation}
(see also \cite[Theorem 6.2]{Di98}).

In the following sections we will mention the results for which
realization functionals were used. We will also present extensions
to weighted spaces and to spaces of multivariate functions.

Like most interesting concepts, realizations were discussed before
the concept was introduced formally.  For instance, the equivalence
\begin{equation}\label{Eq5.12}
R_{r,\varphi  }(f,n^{-r})_p\approx \omega  ^r_\varphi  (f,1/n)_p
\q\q 1\le p\le\infty
\end{equation}
with $P_n$ of (\ref{Eq5.2}) or (\ref{Eq5.3}) was shown already in
\cite{Di-To87} and its trigonometric analogue much earlier.  This
should not diminish the significance of the systematic treatment of
realizations and their importance for various spaces and
applications (not only in relation to algebraic polynomial
approximation).

\section{Sharp Marchaud and sharp converse inequalities}\label{Sec6}

The converse inequality of (\ref{Eq2.6}) is given by
\begin{equation}\label{Eq6.1}
\omega  ^r_\varphi  (f,t)_p\le M(r)t^r \sum^{\floor{\frac 1t}}_{n=1}
n^{r-1}E_n(f)_p, \q 1\le p\le \infty
\end{equation}
with $E_n(f)_p$ given in (\ref{Eq2.7}) was proved in \cite[Theorem
7.2.4]{Di-To87}.  (Note that we write here $n$ instead of $n+1$ in
\cite{Di-To87} as here $\Pi_n =\;\text{\rm span}\,
(1,\dots,x^{n-1}).)$ The Marchaud inequality
\begin{equation}\label{Eq6.2}
\omega  ^r_\varphi  (f,t)_p\le Ct^r \Big\{\int^c_t \;\frac{\omega
^r_\varphi  (f,u)_p}{u^{r+1}}\; du +\Vert  f\Vert  _p\Big\}, \q 1\le
p\le \infty
\end{equation}
was proved in \cite[Theorem 4.3.1]{Di-To87} for a general class of
step weights $\varphi  (x).$ (Not just $\varphi  (x)
=\sqrt{1-x^2}\,.)$

For trigonometric polynomials A.~Zygmund \cite{Zy} and M.~Timan
\cite{TiM58} proved
\begin{equation}\label{Eq6.3}
\omega  ^r(f,t)_{L_p(T)} \le M(r)t^r \Big\{\sum_{n=1}^{\floor{\frac 1t}}
n^{rq-1} E^*_n(f)^q_p\Big\}^{1/q}, \q 1\le p<\infty  , \q q=\min\,(p,2)
\end{equation}
where $\omega  ^r(f,t)_p$ and $E^*_n(f)_p$ are given by (\ref{Eq2.3})
and (\ref{Eq2.4}) respectively.  In addition, it was shown in
\cite{Zy} and \cite{TiM58} that for $1\le p <\infty  $
\begin{equation}\label{Eq6.4}
\omega  ^r(f,t)_{L_p(T)} \le Ct^r \Big[\Big\{\int^c_t\;\frac{\omega
^{r+1}(f,u)^q_{L_p(T)}}{u^{qr+1}} \; du\Big\}^{1/q} +\Vert  f\Vert
_{L_p(T)}\Big] ,  \q q=\;\min\,(p,2).
\end{equation}
(The term $\Vert  f\Vert  _{L_p(T)}$ in (\ref{Eq6.4}) is redundant.)
The classic converse and Marchaud inequalities, i.e. (\ref{Eq6.3})
and (\ref{Eq6.4}) for $1\le p\le \infty  $ when $q=1$
replaces $q=\min(p,2)$,
are clearly
weaker for $1<p<\infty  $ than (\ref{Eq6.3}) and (\ref{Eq6.4}) with
$q=\,\min\, (p,2).$  Moreover, $q=\,\min\, (p,2)$ is the optimal
power in (\ref{Eq6.3}) and (\ref{Eq6.4}) for $1\le p<\infty  .$
Using partially a new proof and extension of (\ref{Eq6.4}) given in
\cite{Di88}, Totik proved in \cite{To88} for $1<p<\infty  $ that
\begin{equation}\label{Eq6.5}
\omega  ^r_\varphi  (f,t)_p\le Ct^r\Big[\Big\{\int^c_t\,
\frac{\omega  ^{r+1}_\varphi  (f,u)^q_p}{u^{rq+1}}\, du\Big\}^{1/q}
+\Vert  f\Vert  _p\Big]\q\text{\rm where} \q q=\,\min\,(p,2)
\end{equation}
(here $\Vert  f\Vert  _p$ can be replaced by $E_{r-1}(f)_p$ but not
eliminated), and he deduced from it for $1<p<\infty  $
\begin{equation}\label{Eq6.6}
\omega  ^r_\varphi  (f,t)_p \le M(r) t^r \Big[\sum^{\floor{\frac 1t}}_{n=1} \,
n^{rq-1}E_n(f)^q_p\Big]^{1/q} \q \text{where} \q q=\,\min\,(p,2).
\end{equation}
In Totik's paper (see \cite{To88}) (\ref{Eq6.5}) is given for $1<
p\le 2$ with a more general step weight $\varphi  .$ For $2<p<\infty
 $ he gave (\ref{Eq6.5}) and (\ref{Eq6.6}) only for $\varphi  (x)
=\sqrt{1-x^2}\,.$

Examples were given in \cite[Section 10]{Da-Di-Ti} to show that the
power $q$ in (\ref{Eq6.5}) and (\ref{Eq6.6}) is optimal for $1<p<4.$
(The power $q$ is probably optimal in (\ref{Eq6.5}) and
(\ref{Eq6.6}) for all $1<p<\infty  .)$

Later it was shown in \cite[Theorem 1.1]{Di-Ji-Le} that
(\ref{Eq6.6}) is valid for $0<p<1$ as well.  Using (\ref{Eq2.6})
which was proved in \cite[Theorem 1.1]{De-Le-Yu} for $0<p<1$ and
applying it to $k+1$ (instead of $k),$ one has (\ref{Eq6.5}) also for
$0<p<1.$

Recently, (see \cite[Theorem 6.2]{Da-Di05}) it was shown that for
$\alpha  <\beta  ,$ $1\le p<\infty  ,$ $q=\,\min\,(p,2)$ and $P(D) =
-\,\frac{d}{dx}\, (1-x^2)\, \frac{d}{dx}$ (among other operators)
one has
\begin{equation}\label{Eq6.7}
K_{2\alpha  }\big(f,P(D)^\alpha  ,t^{2\alpha  })_p\le Ct^{2\alpha
}\Big\{\int^c_t \;\frac{K_{2\beta  }\big(f,P(D)^\beta  ,u^{2\beta
}\big)^q_p}{u^{2\alpha  q+1}} \;du\Big\}^{1/q}
\end{equation}
with $K_{2\alpha  }\big(f,P(D)^\alpha  ,t^{2\alpha  }\big)_p$ given
in (\ref{Eq4.10}).

As we have for $1\le p\le \infty  $ and all $\gamma  >0$
\begin{equation}\label{Eq6.8}
K_{2\gamma  }\big(f,P(D)^\gamma  ,n^{-2\gamma  }\big)_p \le
Cn^{-2\gamma  }\,\sum^n_{k=1} \, k^{2\gamma  -1}E_k(f)_p,
\end{equation}
the inequality (\ref{Eq6.7}) used for $\gamma  =\beta  $ implies
\begin{equation}\label{Eq6.9}
K_{2\alpha  }\big(f,P(D)^\alpha  ,t^{2\alpha  }\big)_p \le
Ct^{2\alpha  }\Big\{\sum_{1\le k\le 1/t} k^{2\alpha  q-1}
E_k(f)^q_p\Big\}^{1/q}.
\end{equation}
For $1<p<\infty  $ and $2\ell =r$ we have
$$
K_{2\ell}\big(f,P(D)^\ell,t^{2\ell}\big)_p \approx \omega
^{2\ell}_\varphi  (f,t)_p
$$
(see \cite[Theorem 7.1]{Da-Di05} and \cite[Chapter 9]{Di-To87}).  In
fact, one can use (\ref{Eq6.8}) and (\ref{Eq6.9}) to obtain
(\ref{Eq6.5}) and (\ref{Eq6.6}).  However, in my opinion, the main
advantage of the technique in \cite{Da-Di05} for polynomial
approximation is not its
applicability to fractional $\alpha  $ but that this method is
applicable to $L_p[-1,1]$ with Jacobi-type weights (see Section~10).

\section{Moduli of smoothness of functions and of their derivatives}\label{Sec7}

For $f,f^{(k)}\in L_p(T),$ $1\le p\le\infty  ,$ it is well-known
(see \cite[p.~46]{De-Lo}) that
\begin{equation}\label{Eq7.1}
\omega  ^r(f,t)_p \le Ct^k \omega
^{r-k}\big(f^{(k)},t\big)_p\q\text{\rm where}\q 1\le k\le r
\end{equation}
and that (see \cite[p.~178]{De-Lo})
\begin{equation}\label{Eq7.2}
\omega  ^{r-k}\big(f^{(k)},t\big)_p \le C\int^t_0\;\frac{\omega
^r(f,u)_p}{u^{k+1}}\;du\q\text{\rm where}\q 1\le k< r.
\end{equation}
It was shown recently (see \cite{Di-Ti}) that the converse-type
inequality (\ref{Eq7.2}) can be improved for $1<p<\infty  $ and has
an analogue for $0<p<1.$

For $\omega  ^r_\varphi
(f,u)_{L_p[-1,1]}$ it was proved in \cite[Theorem 6.2.2 and Theorem
6.3.1]{Di-To87} that
\begin{equation}\label{Eq7.3}
\Omega  ^{r}_\varphi  \big(f,t\big)_p \le Ct^k \omega
^{r-k}_\varphi  \big(f^{(k)},t\big)_{p,\varphi  ^k}\q\text{\rm for}\q
1\le p\le \infty  \q\text{\rm and}\q r>k
\end{equation}
and
\begin{equation}\label{Eq7.4}
\Omega  ^{r-k}_\varphi  \big(f^{(k)},t\big)_{p,\varphi  ^k} \le
C\,\int^t_0\;\frac{\Omega  ^{r}_\varphi  (f,u)_p}{u^{k+1}}\;
du\q\text{\rm for}\q 1\le p\le \infty  \q\text{\rm and} \q r>k
\end{equation}
where
\begin{equation}\label{Eq7.5}
\Omega  ^\ell_\varphi  (g,t)_{p,\varphi  ^m} =\us{\vert  h\vert
<t}\sup\; \Vert  \Delta  ^\ell_{h\varphi  }g\Vert  _{L_{p,\varphi
 ^m}[I(h,\ell)]},
\end{equation}
\begin{equation}\label{Eq7.6}
I(h,\ell) = [-1+2h^2\ell^2,1-2h^2\ell^2]
\end{equation}
and
\begin{equation}\label{Eq7.7}
\Vert  F\Vert  _{L_{p,w}(D)} =\Big\{\int_D \vert  F(x)\vert  ^p
w(x)dx\Big\}^{1/p}
\end{equation}
(and $\Delta  ^\ell_{h\varphi  }f(x)$ is still defined by
(\ref{Eq1.2}) with the underlying interval $[-1,1]).$  In
\cite[Section~5]{Di-Ti}, generalization of (\ref{Eq7.4}) was
achieved, i.e. for $0< p<\infty  $ and $r>k$
\begin{equation}\label{Eq7.8}
\Omega  ^{r-k}_\varphi  \big(f^{(k)},t\big)_{p,\varphi  ^k} \le
C\,\Big\{\,\int^t_0\,\frac{\omega  ^r_\varphi  (f,u)^q_{p}}{u^{qk+1}}\;
du\Big\}^{1/q}
\end{equation}
where $q=\min\,(p,2).$  Simple examples can be given to show that
(\ref{Eq7.3}) does not hold for $0<p<1.$  It can be noted that a
best approximation version of (\ref{Eq7.8}) follows from the proof
in \cite[Section~5]{Di-Ti}, that is,
\begin{equation}\label{Eq7.9}
\Omega  ^{r-k}_\varphi  \big(f^{(k)},t\big)_{p,\varphi  ^k} \le
C\Big\{\sum_{\ell \ge \floor{1/t}} \,\ell^{qk-1} E_\ell(f)^q_p\Big\}^{1/q}
\end{equation}
where $f\in L_p[-1,1],$ $r>k,$ $0<p<\infty  $ and $q=\min\,(p,2).$
In \cite[Section 5]{Di-Ti} it was shown that
$$
\Omega  ^{r-k}_\varphi  \big(f^{(k)},t\big)_{p,\varphi  ^k}\le C
\Big\{\sum_{2^m\ge \floor{1/t}} \,
2^{mkq}E_{2^m}(f)^q_p\Big\}^{1/q}\,,\eqno{(7.9)^\pr}
$$
which is equivalent to (\ref{Eq7.9}).

Another approach to this question which is applicable to Banach
spaces satisfying (\ref{Eq4.8}) (see (\ref{Eq4.5}), (\ref{Eq4.6})
and (\ref{Eq4.7})) is implied by the results in \cite[Sections 6 and
7]{Di98}.  We note that the result below applies to $L_p[-1,1]$ with
$1\le p\le \infty  $ but not with $0< p<1.$

We have for $P(D) =-\,\frac{d}{dx}\, (1-x^2)\;\frac{d}{dx}$ and $B$
satisfying (\ref{Eq4.8})
\begin{equation}\label{Eq7.10}
\begin{aligned}
E_{2\lambda  }(f)_B &\le
\Vert  f-\eta  _\lambda  f\Vert _B\le C E_\lambda  (f)_B,\q
\text{\rm and}\\
E_{2\lambda  } \big( P(D)^\alpha  f\big)_B
&\le \Vert   P(D)^\alpha
f -\eta  _\lambda  \big(P(D)^\alpha  f\big)\Vert  _B \le
C E_\lambda  \big(P(D)^\alpha    f\big)_B
\end{aligned}
\end{equation}
where $\eta  _\lambda  $ is the
 de~la~Vall\'ee Poussin-type operator defined
by (\ref{Eq5.4})
using (\ref{Eq4.6}) and  (\ref{Eq4.7}).  (Other de~la~Vall\'ee
Poussin-type operators will yield a result similar to (\ref{Eq7.10}).)

Using the realization theorem (see \cite[Theorem
7.1]{Di98}) given by
\begin{equation}\label{Eq7.11}
K_{2\alpha  }
\Big(f,\big(-\,\frac{d}{dx}\,(1-x^2)\,\frac{d}{dx}\big)^\alpha
,\frac{1}{n^{2\alpha  }}\Big)_B\approx \Vert  f-V_nf\Vert  _B
+ \frac{1}{n^{2\alpha  }} \;\Big\Vert  \Big(-\,\frac{d}{dx}\,
(1-x^2)\,\frac{d}{dx}\Big)^\alpha  V_nf\Big\Vert  _B,
\end{equation}
with $V_nf = \eta  _{1/n}  f$ or other de~la~Vall\'ee Poussin-type
operators,
one has the following result:

For $\alpha  <\beta  $ and
$\big(-\,\frac{d}{dx}\,(1-x^2)\,\frac{d}{dx}\big)^\alpha  f\in B$ we
have, using
\cite[Theorem 7, (7.10) and (5.11)]{Di98},
\begin{equation}\label{Eq7.12}
\begin{aligned}
K_{2\beta  }
&\Big(f,\big(-\,\frac{d}{dx}\,(1-x^2)\,\frac{d}{dx}\big)^\beta
,t^{2\beta  }\Big)_B  \\
&\le Ct^{2(\beta  -\alpha  )} K_{2(\beta  -\alpha  )}
\,\Big(\big(-\,\frac{d}{dx} \,(1-x^2)\,\frac{d}{dx}\big)^\alpha
f,\big(-\,\frac{d}{dx}\,(1-x^2)\,\frac{d}{dx}\big)^{\beta  -\alpha
},t^{2(\beta  -\alpha  )}\Big)_B.
\end{aligned}
\end{equation}

For $\alpha  <\beta  $ and $P(D) =-\,\frac{d}{dx}\,
(1-x^2)\,\frac{d}{dx}$ we also have, using \cite[Theorem 7.1]{Di98},
\begin{equation}\label{Eq7.13}
K_{2(\beta  -\alpha  )} \Big(P(D)^\alpha  f,\big(P(D)\big)^{\beta
-\alpha  },t^{2(\beta  -\alpha  )}\Big)_B
\le C\int^{t}_0\;\frac{K_\beta  \big(f,P(D)^\beta  ,u^{2\beta
}\big)_B}{u^{2\alpha  +1}}\;du.
\end{equation}

For $B=L_p[-1,1],$ $1<p<\infty  $ one can follow \cite{Da-Di05} and
obtain a sharper version of (\ref{Eq7.13}), that is
\begin{equation}\label{Eq7.14}
\begin{aligned}
K_{2(\beta  -\alpha  )}& \big(P(D)^\alpha  f,P(D)^{\beta
-\alpha  },t^{2(\beta  -\alpha  )}\big)_{L_p[-1,1]}\\
&\le C\Big\{\int^{t}_0\;\frac{K_\beta  \big(f,P(D)^\beta  ,u^{2\beta
}\big)^q_{L_p[-1,1]}}{u^{2\alpha q+1}}\Big\}^{1/q}, \q
q=\min\,(p,2).
\end{aligned}
\end{equation}

For the rate of best approximation $E_n(f)_B$ given in (\ref{Eq2.7})
or (\ref{Eq4.9}) (when (\ref{Eq4.8}) is satisfied),
(\ref{Eq7.13}) and (\ref{Eq7.14}}) take the forms
\begin{equation}\label{Eq7.15}
K_{2(\beta  -\alpha  )} \Big(P(D)^\alpha  f,\big(P(D)\big)^{\beta
-\alpha  },t^{2(\beta  -\alpha  )}\Big)_B
\le C\sum_{\ell\ge \floor{1/t}} \ell^{2\alpha  -1} E_\ell (f)_B,
\end{equation}
and for $1<p<\infty  $
\begin{equation}\label{Eq7.16}
\begin{aligned}
K_{2(\beta  -\alpha  )} \Big(P(D)^\alpha  f,& P(D)^{\beta  -\alpha
},t^{2(\beta  -\alpha  )}\Big)_{L_p[-1,1]}\\
&\le C\Big\{\sum_{\ell \ge \floor{1/t}} \ell^{2\alpha  q-1} E_\ell
(f)^q_{L_p[-1,1]}\Big\}^{1/q}\q \text{\rm where}\q
q = \min\, (p,2)
\end{aligned}
\end{equation}
respectively.

\section{Relations with Bernstein polynomial approximation and
other linear\\ operators}\label{Sec8}

Chapters 9 and 10 of \cite{Di-To87} were dedicated to relations
between $\omega  ^r_\varphi  (f,t)_p$ (with appropriate $\varphi  $
and domain) and the rate of convergence of Bernstein, Szasz and
Baskakov operators (including appropriate combinations and
modifications).

We remind the reader that the Bernstein operator is given by
\begin{equation}\label{Eq8.1}
B_n(f,x) =\sum^n_{k=0} \binom{n}{k} x^k(1-x)^{n-k} f\Big(\frac
kn\Big) \equiv \sum^n_{k=0} P_{n,k} (x)f\Big(\frac
kn\big)\q\text{\rm for} \q x\in [0,1].
\end{equation}

Perhaps the first real progress in the last twenty years was the
general group of concepts called strong converse inequalities S.C.I.
(see \cite{Di-Iv}).  In \cite[Section 8]{Di-Iv} it was shown as
one of the applications of the general method given in \cite[Section
3]{Di-Iv} that
\begin{equation}\label{Eq8.2}
\omega  ^2_\varphi  (f,n^{-1/2})_{C[0,1]} \le
C\Big(\Vert  B_n f-f\Vert  _{C[0,1]} + \Vert  B_{An} f-f\Vert
_{C[0,1]}\Big)
\end{equation}
for some $A>1$ where $\omega  ^2_\varphi  (f,t)_{C[0,1]}$ (in
relation to Bernstein polynomials) is a copy of $\omega  ^2_\varphi
(f,t)_{C[-1,1]}$ given by (\ref{Eq1.1}) and (\ref{Eq1.2}) in which
$[0,1]$ replaces $[-1,1]$ and $\varphi  (x) =\sqrt{x(1-x)}$ replaces
$\sqrt{1-x^2}\,.$  The inequality (\ref{Eq8.2}) is a strong converse
inequality of type $B$ (with two terms on the right hand side) and
is called ``strong'' as it matches the direct result (see
\cite{Di-To87}) given by
\begin{equation}\label{Eq8.3}
\Vert  B_nf - f\Vert  _{C[0,1]} \le C\omega  ^2_\varphi
\Big(f, n^{-1/2}\Big)_{C[0,1]}.
\end{equation}
One observes that (\ref{Eq8.2}) implies
$$
\omega  ^2_\varphi  (f,n^{-1/2})_{C[0,1]} \le C\,\us{k\ge
n}\sup\;\Vert  B_kf-f\Vert  _{C[0,1]}, \eqno{(8.2)^\pr}
$$
which is a strong converse inequality of type $D$ in the terminology
of \cite{Di-Iv}.  Combining $(8.2)^\pr$ with (\ref{Eq8.3}), one has
$$
\omega  ^2_\varphi  (f,n^{-1/2})_{C[0,1]} \approx \us{k\ge n}\sup\;
\Vert  B_nf -f \Vert  _{C[0,1]}\,.
$$

In \cite[Remark 8.6]{Di-Iv} it was conjectured that the superior
strong converse inequality of type A is also valid, that is, that
\begin{equation}\label{Eq8.4}
\omega  ^2_\varphi  (f,n^{-1/2})_{C[0,1]} \le
C\Vert  B_n f-f\Vert  _{C[0,1]}
\end{equation}
which, together with (\ref{Eq8.3}), implies
\begin{equation}\label{Eq8.5}
\Vert  B_n f-f\Vert  _{C[0,1]} \approx \omega  ^2_\varphi
\Big(f,n^{-1/2}\Big)_{C[0,1]}.
\end{equation}

In a remarkable paper (see \cite{To94}) V. Totik gave the first proof of
(\ref{Eq8.4}). He used an intricate modification of the parabola
technique.  Totik's method is applicable to Bernstein, Szasz and
Baskakov operators.  Explicitly, Totik treated the Szasz-Mirakian
operator given by
\begin{equation}\label{Eq8.6}
S_n (f,x) =\sum^\infty  _{k=0} e^{-nx} \,\frac{(nx)^k}{k!}\,
f\big(\frac kn\big),
\end{equation}
for which he showed
\begin{equation}\label{Eq8.7}
\Vert  S_n(f,x) -f(x)\Vert  _{C[0,\infty  )} \approx \omega
^2_\varphi  \Big(f, n^{-1/2}\Big)_{C[0,\infty  )}
\end{equation}
where $\omega  ^2_\varphi  (f,t)_{C[0,\infty  )}$ is defined on
$[0,\infty  )$ (instead of $[-1,1])$ and $\varphi  (x) = \sqrt x$
(instead of $\sqrt{x(1-x)}$ or $\sqrt{1-x^2}\,).$  The proof of
(\ref{Eq8.7}) is neater than that of (\ref{Eq8.4}) as $[0,\infty  )$
has only one finite endpoint and $\sqrt x$ is simpler than
$\sqrt{x(1-x)}\,.$  Totik stated that the proof in the case of
Bernstein and Baskakov operators is essentially the same.  To prove
(\ref{Eq8.4}) directly would be just a bit longer, more cluttered
and would perhaps obscure the idea.

The second proof of (\ref{Eq8.4}) was given by Knopp and Zhou (see
\cite{Kn-Zh94}), who used the fact that
\begin{equation}\label{Eq8.8}
\frac 1n\;\Big\Vert  \varphi  ^2\Big(\frac{d}{dx}\Big)^2 B^m_n
f\Big\Vert  _{C[0,1]} \le C(m)\Vert  f\Vert  _{C[0,1]}
\end{equation}
with $C(m)$ small enough for some $m$ (independent of $n$ and $f)$
being sufficient. $(B^m_n f =
B_n B^{m-1}_nf$ and $B^1_n f = B_n f.)$
It was shown in \cite[Section 4]{Di-Iv} for a large class of
operators $O_n$ and an appropriate differential operator $P(D)$ that a condition like
$\Vert  P(D)O^m_n f\Vert  _B \le C(m)\Vert  f\Vert  _B$ would be
sufficient for proving S.C.I. of type A provided that $C(m)$ is
small enough.  In \cite{Kn-Zh95} (which precedes \cite{Kn-Zh94}) a
general ingenious method was given to show that under some
conditions $C(m)\to 0$ as $m\to\infty  $ for many operators. This
technique is useful and the conditions necessary are easy to verify
when the various operators treated commute, and it is applicable to
many spaces (not just $L_\infty  ).$  However, as $B_n B_m f\ne B_m
B_n f$ and $\varphi  ^2\big(\frac{d}{dx}\big)^2 B_n (f,x)\ne B_n
(\varphi  ^2 f^{\pr\pr}, x)$ even for very smooth functions, the
proof in \cite{Kn-Zh94} becomes extremely complicated.  I note that
in papers of X.~Zhou with Knoop and others strong converse
inequalities are called lower estimate (to match the direct result
like (\ref{Eq8.3}) which Zhou et~al. call the upper estimate).
Besides this linguistic innovation, and their new idea to show $C(m)
=o(1)$ as $m\to\infty  ,$ they also repeated  the
arguments of \cite[Sections 3-4]{Di-Iv}, perhaps because they felt
they could explain things better.

The third proof of (\ref{Eq8.4}), given by C.~Sanguesa (see \cite{Sa}), uses
probabilistic ideas to show that $C(m)$ of (\ref{Eq8.8}) is
sufficiently small for $m=3.$  The ideas of \cite{Sa} can be
translated from probabilistic to classical analytic.

While S.C.I. of type B are now quite easy to prove and yield most
results about the relation between the $K\text{\rm -functional}$ and
$\Vert  O_nf-f\Vert  ,$ S.C.I. of type A are much more elegant and
hence more desirable.  (They are also more amenable to iterations.)
 I still would like to see a new simple proof of (\ref{Eq8.4}) which
I am sure will have implications for other operators. One wonders
what condition on the sequence of operators (not just the Bernstein
polynomials), which is easy to verify, is sufficient to guarantee
that a S.C.I. of type B implies a S.C.I. of type A.

As the Bernstein operators are not defined on $L_p[0,1]$ for $1\le
p<\infty  ,$ their Kantorovich modification given by
\begin{equation}\label{Eq8.9}
K_n(f,x) =\sum^n_{k=0} \,\binom{n}{k} \,x^k(1-x)^{n-k}
\Big[(n+1)\int^{(k+1)/(n+1)}_{k/(n+1)} f(u)du\Big]
\end{equation}
was extensively used. (Similar extensions were given to Szasz and Baskakov
operators.)

In \cite{Go-Zh} the following S.C.I. of type A is claimed for $1\le
p\le \infty  :$
\begin{equation}\label{Eq8.10}
\Vert  K_n f-f\Vert  _{L_p[0,1]} \approx \inf \Big(\Vert  f-g\Vert
_{L_p[0,1]} +\frac 1n\,\Big\Vert  \frac{d}{dx}\,
x(1-x)\,\frac{d}{dx}\, g\Big\Vert  _{L_p[0,1]}\Big).
\end{equation}

One recalls that the affine transformation $[-1,1]\to [0,1]$ and
(\ref{Eq4.2}), (\ref{Eq4.3}) and
(\ref{Eq4.4}) here imply for $1<p<\infty  $
\begin{equation}\label{Eq8.11}
\omega  ^2_\varphi  \Big(f,\,\frac{1}{\sqrt n}\Big)_{L_p[0,1]}
+\frac 1n\;\Vert  f\Vert  _{L_p[0,1]} \approx \inf \,\Big(\Vert
f-g\Vert  _{L_p[0,1]} +\frac 1n\, \Big\Vert  \frac{d}{dx}\,
x(1-x)\;\frac {d}{dx}\,g\Big\Vert   _{L_p[0,1]} \Big).
\end{equation}
For $p=1$ and $p=\infty  $ (\ref{Eq8.11}) is not valid (see
\cite[p.~88]{Da-Di05}).

Most of the (multitude of) papers on Bernstein-type operators deal
with:

\vs
\begin{description}
\item  {(a)} Combinations (for higher levels of smoothness).

\item  {(b)} Weighted approximation of the operators (see also
Sections
\ref{Sec10} and \ref{Sec14}).

\item  {(c)} Different step-weights (see also Section \ref{Sec14}).

\item  {(d)} Multivariate analogues (see also Section \ref{Sec12}).

\item  {(e)} Simultaneous approximation.

\item {(f)} Shape-preserving properties (see also Section
\ref{Sec15}).

\item  {(g)} Other modifications and generalizations.
\end{description}

If I describe all related results on the subject, I will exhaust
both myself and the reader (who is probably tired already), and
therefore I will try to be somewhat more selective in this survey.
Even after remarks in the following sections, the treatment is by no
means complete and many, perhaps most, results on the topics
(a) -- (g) are not described.

The Bernstein polynomial operator preserves many properties. Its
rate of convergence is equivalent to $\omega  ^2_\varphi
\big(f, n^{-1/2}\big)_{C[0,1]}.$  Realization results using
it are valid (and weaker than (\ref{Eq8.5})).  Moreover, the
Bernstein polynomial operator is a model for many other operators,
mostly yielding similar or weaker results for $C[0,1].$  Therefore,
it was a surprise that a modification emerged that had many ``nice''
properties, some different from those of $B_nf,$ yet extremely
useful.  Such an operator, introduced by Durrmeyer (see \cite{Du}
and \cite{De81}), is now called the Durrmeyer-Bernstein polynomial
operator and is given by
\begin{equation}\label{Eq8.12}
M_n(f,x) =\sum^n_{k=0} P_{n,k}(x)(n+1)\int^1_0 P_{n,k}(y)f(y)dy, \q
P_{n,k}(x) =\binom{n}{k}\,x^k(1-x)^{n-k}.
\end{equation}

Among the properties of $M_n(f,x)$ we state:

\begin{description}
\item  {I.} $M_nf = M_n(f,x): L_p[0,1]\to \Pi_{n+1}$ for $1\le p\le
\infty  .$
\item  {II.} $\Vert  M_nf\Vert  _{L_p[0,1]} \le \Vert  f\Vert
_{L_p[0,1]}$ for $1\le p\le \infty  .$
\item  {III.} $\la M_n f,g\ra = \la f,M_n g\ra$ where $\la F,G\ra
=\int^1_0 F(x)G(x)dx.$
\item  {IV.} For $f\sim \overset\infty  {\us{k=0}\sum} P_k f,\q M_nf\sim
\overset n{\us{k=0}\sum} a_k P_k f$
 where $P_kf$ is given by (\ref{Eq4.6}) with ${\cal D}=[0,1]$ and
(\ref{Eq4.7}) is replaced by
\end{description}
\begin{equation}\label{Eq8.13}
\frac{d}{dx}\,x(1-x)\,\frac{d}{dx}\,\varphi  _k(x) = -k(k+1)\varphi
_k(x), \q \int^1_0 \varphi  _k(x)\varphi  _\ell (x)dx =
\begin{cases} 0 &k\ne\ell,\\ 1 & k=\ell.\end{cases}
\end{equation}

\noi
As a result of IV one has:

\begin{description}
\item  {V.} $M_nM_kf = M_k M_n f.$
\item  {VI.} $\frac{d}{dx}\,\big(x(1-x)\big)\,\frac{d}{dx}\,M_nf =
M_n\big(\frac{d}{dx}\,x(1-x)\,\frac{d}{dx}\,f\big)$ for $f$ smooth
enough.
\item  {VII.} $M_n f - f =\overset\infty  {\us{k=n+1}\sum}
\,\frac{1}{k(k+1)}\,\frac{d}{dx}\,\big(x(1-x)\big)\,\frac{d}{dx}\,
M_kf.$
\end{description}

Using all these properties, it was shown in \cite[Theorem
6.3]{Ch-Di-Iv} for $1\le p\le \infty  $ that
\begin{equation}\label{Eq8.14}
\Vert  M_n f-f\Vert  _p\approx \inf \Big(\Vert  f-g\Vert  _p +\frac
1n\,\Big\Vert  \frac{d}{dx}\,
\big(x(1-x)\big)\,\frac{d}{dx}\,g\Big\Vert  _p\Big).
\end{equation}

Many properties of $M_nf$ were investigated and the proof did not
always use the obvious advantages enumerated above (by I $\to$ VII).

Other multiplier-type polynomial approximation processes are the
Ces\`aro means $C^\ell_n(f,x)$ given in (\ref{Eq4.5}) with
$P_kf = P_k(f,x)$ given in (\ref{Eq4.6}) and
$\varphi
_k(x)$ given above (in (\ref{Eq8.11})) and the Riesz means
\begin{equation}\label{Eq8.15}
R_{\lambda  ,\alpha  ,\ell}f =\sum_{\lambda  (k)<\lambda  }
\Big(1-\big(\frac{\lambda  (k)}{\lambda  }\big)^\alpha  \Big)^\ell
P_k f, \q \lambda  (k) =k(k+1).
\end{equation}

F. Dai proved in \cite{Da} for $\ell\ge 1$ and $1\le p\le\infty  $
that
\begin{equation}\label{Eq8.16}
\Vert  C^\ell_n f-f\Vert  _{L_p[0,1]} \approx \inf \,\Big(\Vert
f-g\Vert  _{L_p[0,1]} +\frac 1n\, \big\Vert
\big(P(D)\big)^{1/2}g\big\Vert_{L_p[0,1]}\Big)
\end{equation}
and
\begin{equation}\label{Eq8.17}
\Vert  R_{n^2,\alpha  ,\ell} f-f\Vert  _{L_p[0,1]} \approx
\inf\,\Big(\Vert  f-g\Vert  _{L_p[0,1]} +\frac {1}{n^{2\alpha  }}\,
\big\Vert  \big(P(D)\big)^\alpha  g\big\Vert \Big) _{L_p[0,1]}
\end{equation}
with $P(D) = -\,\frac{d}{dx}\,\big(x(1-x)\big)\,\frac{d}{dx}\,.$

We note that in this section we use linear operators and S.C.I.
which, when applicable, are more powerful than results on
$K\text{\rm-functionals}$ or realizations.

\section{Weighted moduli of smoothness, doubling weights}\label{Sec9}

In a series of articles Mastroianni and Totik introduced the concept
of doubling weights and showed that many results about trigonometric
polynomials on $T$ and about algebraic polynomials on $[-1,1]$ can
be extended (or modified) to include weighted $L_p$ versions with such
weights.  I will deal here only with results for algebraic polynomials
on $[-1,1].$

Mastroianni and Totik also gave results related to
earlier concepts such as the Muckenhoupt $A_p$ condition and
others.  Let me now briefly describe the concepts involved.

A doubling weight on $[-1,1]$ is a non-negative measurable function
$w(x)$ satisfying
\begin{equation}\label{Eq9.1}
w(2I) \equiv \int_{2I\cap[-1,1]} w(t)dt \le L\int_I w(t)dt \equiv
Lw(I)
\end{equation}
where $I\subset [-1,1],$ $2I$ is the interval with the same
midpoint and twice the length of $I,$ and $L$ is the doubling
constant.  In \cite[Lemma 2.1]{Ma-To00} many definitions equivalent
to (\ref{Eq9.1}) were given.

A non-negative measurable function $w(x)$ is a weight satisfying the
$A_\infty  $ condition if for any set $E,$ $E\subset I\subset
[-1,1]$ with $m(E)\equiv \vert  E\vert  \ge \alpha  \vert  I\vert  $
\begin{equation}\label{Eq9.2}
w(E) \equiv \int_E w(t)dt\ge \beta  w(I)
\end{equation}
with $\beta  =\beta  (\alpha  ).$

A non-negative measurable weight function $w(x)$ satisfies the $A_p$
condition, for some $p,$ $1\le p<\infty  ,$ if for $q= p/(p-1)$
\begin{equation}\label{Eq9.3}
\Big(\frac{1}{\vert  I\vert  }\,\int_I
w(t)dt\Big)\,\big(\frac{1}{\vert  I\vert  }\, \int_I w(t)^{-q/p}
dt\Big)^{p/q} \le A
\end{equation}
for all $I\subset [-1,1].$

A non-negative measurable weight function $w(x)$ satisfies the $A^*$
condition if
\begin{equation}\label{Eq9.4}
w(x)\le L\;\frac{1}{\vert  I\vert  }\, \int_I w(t)dt
\end{equation}
for $x\in I\subset [-1,1]$ and $L$ independent of $x$ and $I.$
(Note that satisfying the $A^*$ condition implies that $w(x)$ is
bounded.)

Clearly, the conditions are ordered in increasing strength and the
doubling weight condition is the most general (weakest).

In the next section we will describe results for the Jacobi weights
which are not known, not valid, or just not applicable to the
classes of weights mentioned above.  We note that except for $A^*$
the above-mentioned classes of weights contain the Jacobi weights
treated
in
Section \ref{Sec10}, and $A^*$ contains the bounded Jacobi weights.
Furthermore, we note that, for example, the weight
$$
w(t) = h(t)\,\prod^k_{j=1}\,\vert  t-x_j\vert  ^{\gamma  _j},
$$
with $\gamma  _j >-1,$ $x_j\in [-1,1],$ $x_j <x_{j+1}$  and $h(t)$ a positive measurable
function satisfying $0 <A\le h(t)\le B<\infty  ,$ is a doubling
weight.

One defines $w_n(x)$ by
\begin{equation}\label{Eq9.5}
w_n(x) =\frac{1}{\Delta  _n(x)}\,\int_{\big(x-\Delta  _n(x),x+\Delta
 _n(x)\big)\cap[-1,1]} w(u)du, \q \Delta  _n(x) =\frac{\sqrt{1-x^2}}{n} +
\frac{1}{n^2}
\end{equation}
and notes that $w_n(x)$ is a doubling weight whenever $w(x)$ is.

We denote (as usual)
\begin{equation}\label{Eq9.6}
\Vert  f\Vert  _{L_p(w)} = \Vert  f\Vert  _{w,p} =\Big\{\int^1_{-1}
\vert  f(x)\vert  ^p w(x)dx\Big\}^{1/p}, \q 0<p<\infty  ,
\end{equation}
$$
\Vert  f\Vert  _{L_\infty  (w)} = \Vert  f\Vert  _{w,\infty  } =
\;\us{x\in [-1,1]} {\text{\rm ess sup}}\,\vert  f(x)w(x)\vert
,\eqno(9.6)^\pr
$$
\begin{equation}\label{Eq9.7}
E_n(f)_{w,p} \equiv \us{P_n\in \Pi_n}{\text{\rm inf}}\, \Vert
f-P_n\Vert  _{w,p}, \q \Pi_n = \;\text{\rm span}\, (1,\dots,x^{n-1})\, ,
\end{equation}
and
\begin{equation}\label{Eq9.8}
\omega  ^r_\varphi  (f,t)_{w,p} =\us{\vert  h\vert  \le t}{\text{\rm
sup}}\, \Vert  \Delta  ^r_{h\varphi  }f\Vert  _{w,p}
\end{equation}
where $\Delta  ^r_{h\varphi  }f$ is given in (\ref{Eq1.2}).

A Jackson-type result for general doubling weights and $1\le
p<\infty  $ was given
(see \cite[Theorem 3.2]{Ma-To98}) by
\begin{equation}\label{Eq9.9}
E_n(f)_{w,p}\le \frac{C}{n^r} \,\Vert  f^{(r)}\varphi  ^r_n\Vert
_{w_n,p}\,, \q \varphi  _n(x) = \sqrt{1-x^2} +\frac 1n
\end{equation}
where $w_n$ is given in (\ref{Eq9.5}) and $f,\dots,f^{(r-1)}\in
\;\text{\rm A.C.}_{\ell\text{\rm oc}}.$  For a weight $w$ satisfying
the $A_p$ condition $w$ can replace $w_n$ in (\ref{Eq9.9}) (see
\cite[Theorem 3.4]{Ma-To98}), and when $w(x)\approx w_n(x)$ for
$x\in [-1+\frac{1}{n^2}, 1-\frac{1}{n^2}],$ both $w_n$ and $\varphi
_n$ can be replaced by $w$ and $\varphi  $ in (\ref{Eq9.9}) (see
\cite[Theorem 3.6]{Ma-To98}).  For $p=\infty  $ a Jackson-type
result was given by
$$
E_n(f)_{w_n,\infty  } \le \frac{C}{n^r} \, \Vert  f^{(r)}\varphi
^r\Vert  _{w_n,\infty  } \eqno(9.9)^\pr
$$
(see \cite[Theorem 1.1]{Ma-To99}).

Clearly, (\ref{Eq9.9}) and $(9.9)^\pr$ imply
\begin{equation}\label{Eq9.10}
E_n(f)_{w,p} \le C\;\us g{\text{\rm inf}}\, (\Vert  f-g\Vert
_{w_n,p} + n^{-r}\Vert  g^{(r)} \varphi  ^r_n\Vert  _{w_n,p})
\equiv CK_{r,\varphi  _n}(f,n^{-r})_{w_n,p}
\end{equation}
for $1\le p<\infty  ,$ and
$$
E_n(f)_{w_n,\infty  } \le C\;\us g{\text{\rm inf}}\, (\Vert  f-g\Vert
_{w_n,\infty  } + n^{-r}\Vert  g^{(r)} \varphi  ^r\Vert  _{w_n,\infty
 }) \equiv CK_{n,\varphi  }(f,n^{-r})_{w_n,\infty  }.
\eqno(9.10)^\pr
$$

We note that the price for dealing with such general weights as the
doubling weight is that in (\ref{Eq9.10}) and $(9.10)^\pr$ we do
not have one $K\text{\rm -functional}$ but a sequence of (somewhat)
different ones which depend on $n.$  (Recall that for $w=1,$
(\ref{Eq2.6}) and (\ref{Eq3.2}) imply $E_n(f)_p \le CK_{r,\varphi
}(f,n^{-r})_p$ for $1\le p \le \infty  .)$  Using \cite[Theorem
3.6]{Ma-To98} and \cite[Theorem 1.2]{Ma-To99}, we also have one
$K\text{\rm -functional}$ when $w_n(x)\approx w(x)$for $x\in
[-1+\frac{1}{n^2}, 1-\frac{1}{n^2}]$
and for the class of weights given by $A^*$
as in these cases $w_n$ and $\varphi  _n$
are replaced by $w$ and $\varphi  $ in (\ref{Eq9.9}).

Following the proof in \cite[Theorem 2.1.1]{Di-To87}, one has
$$
\omega  ^r_\varphi  (f,t)_{w_n,\infty  } \approx K_{r,\varphi
}(f,t^r)_{w_n,\infty  } \q\text{\rm and} \q \omega  ^r_{\varphi _n}
(f,t)_{w_n,p} \approx K_{r,\varphi  _n}(f,t^r)_{w_n,p}
$$
(see \cite[p.~188]{Ma-To01}).

The converse result
\begin{equation}\label{Eq9.11}
\omega  ^{r+2}_\varphi  \Big(f,\frac 1n\Big)_{w_n,\infty  } \le
Cn^{-r}\sum^n_{k=1} k^{r-1}E_k(f)_{w_k,\infty  }
\end{equation}
was proved in \cite[(1.8)]{Ma-To01} where it was shown that in
general $r+2$ on the left of (\ref{Eq9.11}) cannot be improved. For
$w$ satisfying the $A^*$ condition, $\omega  ^r_\varphi
(f,1/n)_{w_n,\infty  }$ can replace $\omega  ^{r+2}_\varphi
(f,1/n)_{w_n,\infty  }$ in (\ref{Eq9.11}).

Some questions such as: estimating $E_n(f)_{w_n,p}$ by $\omega
^r_{\varphi  _n}(f,t)_{w_n,p}$ for $0<p<1,$ the connection between
$\omega  ^r_{\varphi  _n}(f,t)_{w_n,p}$ and appropriate
realizations, and whether $r+2$ on the left of (\ref{Eq9.11}) is
still necessary for $1\le p <\infty  ,$ were not considered as far
as I know.

A wealth of results about inequalities concerning polynomials on
$[-1,1]$ in weighted $L_p$ norms were given in the series of papers
mentioned and in particular in \cite{Ma-To00}.  These inequalities
will be crucial for further investigations.

For a doubling weight $w$ and for $w_n(x)$ given by (\ref{Eq9.5})
it was shown \cite[Theorem 7.2]{Ma-To00} that for $P_n\in \Pi_n$ and
$1\le p<\infty  $ one has
\begin{equation}\label{Eq9.12}
\frac 1C \, \int^1_{-1} \,\vert  P_n\vert  ^p w\le \int^1_{-1}\,
\vert  P_n\vert  ^p w_n \le C\,\int^1_{-1} \,\vert  P_n\vert  ^pw.
\end{equation}
The Bernstein and Markov inequalities  for a doubling weight (see \cite[Theorem
7.3 and 7.4]{Ma-To00}) were given for $1\le p<\infty  $ and $P_n\in \Pi_n$ by
\begin{equation}\label{Eqnew}
\int^1_{-1} \varphi  ^p \vert  P^\pr_n\vert  ^p w\le C
n^p\int^1_{-1}\vert  P_n\vert  ^p w
\end{equation}
and
\begin{equation}\label{Eq9.13}
\int^1_{-1} \, \vert  P^\pr_n\vert  ^pw\le C
n^{2p}\int^1_{-1}\,\vert  P_n\vert  ^p w, \q 1\le p<\infty
\end{equation}
respectively.

The Nikol'skii inequality was given in two different forms for $1\le
p<q<\infty  $ and $P_n\in \Pi_n$ (see \cite[p.~67]{Ma-To00}) by
\begin{equation}\label{Eq9.14}
\Big(\int^1_{-1} \vert  P_n\vert  ^qw\Big)^{1/q} \le Cn^{\frac 2p
-\frac 2q} \Big(\int^1_{-1} \vert  P_n\vert  ^p w^{p/q}\Big)^{1/p}
\end{equation}
and by
\begin{equation}\label{Eq9.15}
\Big(\int^1_{-1} \vert  P_n\vert  ^q w\Big)^{1/q} \le
Cn^{\frac1p-\frac 1q} \Big(\int^1_{-1} \vert  P_n\vert  ^p w^{p/q}
\varphi  ^{\frac pq-1}\Big)^{1/p}.
\end{equation}
Note that for the special case of Jacobi weights a third different
form will be presented in the next section, and while (\ref{Eq9.14})
and (\ref{Eq9.15}) are best possible of their type, the third form
(for Jacobi weights) will also be best possible. I find it amusing
to see three different Nikol'skii-type inequalities for algebraic
polynomials, all best possible in their way, which treat the weight
on the right hand side differently.

For $w$ satisfying the $A^*$ condition one has (\cite[p.
69]{Ma-To00}) the Bernstein inequality
\begin{equation}\label{Eq9.16}
\Vert  \varphi  P^\pr_n w\Vert  _{L_\infty  [-1,1]} \le Cn\Vert  P_n
w\Vert  _{L_\infty  [-1,1]},
\end{equation}
the Markov-Bernstein inequality
\begin{equation}\label{Eq9.17}
\Vert  P^\pr_n w\Vert  _{L_\infty  [-1,1]} \le Cn^2 \Vert  P_n
w\Vert  _{L_\infty  [-1,1]} ,
\end{equation}
and the Nikol'skii-type inequality
\begin{equation}\label{Eq9.18}
\Vert  P_n w\Vert  _{L_\infty  [-1,1]} \le Cn^{2/p} \Vert  P_n
w\Vert  _{L_p[-1,1]}, \q p <\infty  .
\end{equation}
For Jacobi-type weights (\ref{Eq9.18}) is improved on in the next
section, but as applicable to all $w$ satisfying the $A^*$
condition, the inequality (\ref{Eq9.18}) is best possible as well.

The multivariate analogues were not considered for algebraic
polynomials. (For the multivariate situation on the sphere results
using spherical harmonic polynomials are treated in \cite{Da06}.)
The case $0<p<1$ for algebraic polynomials was not considered
explicitly.  For trigonometric polynomials analogues of Bernstein,
Marcinkiewicz, Nikol'skii and Schur type inequalities (but not the
Jackson-type inequality) are given for $0<p$ and doubling weights or
$A^*$ weights in \cite{Er}.  These results can probably be extended
to algebraic polynomials with appropriate modifications but Erd\'elyi
states ``For technical reasons we discuss only the trigonometric
cases''.

\section{Weighted moduli, Jacobi-type weights}\label{Sec10}

The Jacobi weights  given by
\begin{equation}\label{Eq10.1}
w(x) = w_{\alpha  ,\beta  }(x) = (1-x)^\alpha  (1+x)^\beta  , \q
\alpha  >-1, \q \beta  >-1,
\end{equation}
are doubling weights, (that is, they satisfy (\ref{Eq9.1})), and for
$\alpha  \ge 0,$ $\beta  \ge 0$ they are also $A^*$ type weights
(i.e. satisfying (\ref{Eq9.4})).  Moreover, for $x\in
[-1+\frac{1}{n^2}, 1-\frac{1}{n^2}]$ they satisfy $w(x)\approx
w_n(x),$ and hence the discussion in the last section implies for
$1\le p<\infty  $ and $w(x)$ (see \cite[Theorem 3.6]{Ma-To98})
\begin{equation}\label{Eq10.2}
\begin{aligned}
E_n(f)_{w,p} &\le C\;\us g{\text{inf}} \,\big(\Vert  f-g\Vert  _{w,p} +
n^{-r}\Vert  g^{(r)}\varphi  ^{r}\Vert  _{w,p}\big)\\
&\equiv CK_{r,\varphi  }(f,n^{-r})_{w,p}.
\end{aligned}
\end{equation}
For $\alpha  \ge 0,$ $\beta  \ge 0$ (\ref{Eq10.2}) follows for
$p=\infty  $ as well (see \cite[Theorem 1.2]{Ma-To99}.

For the Jacobi weights different $K\text{\rm -functionals}$
(see \cite{Ch-Di97}, \cite{Di98} and
\cite{Da-Di05}), which are given for $\alpha  >-1,$ $\beta  >-1$ by
\begin{equation}\label{Eq10.3}
K_\gamma  \big(f,P_{\alpha  ,\beta  }(D)^\gamma  ,t^{2\gamma
}\big)_{w_{\alpha  ,\beta  },p} = \;\text{\rm inf}\, \big(\Vert
f-g\Vert  _{w_{\alpha  ,\beta  },p} + t^{2\gamma  }\Vert  P_{\alpha
,\beta  }(D)^\gamma   g\Vert  _{w_{\alpha  ,\beta  },p}\big)
\end{equation}
where
\begin{equation}\label{Eq10.4}
P_{\alpha  ,\beta  }(D) = -w_{\alpha  ,\beta  }(x)^{-1} \,
\frac{d}{dx}\; w_{\alpha  ,\beta  }(x) (1-x^2)\;\frac{d}{dx}
\end{equation}
were shown to be useful.

We note that for integer $\gamma  $ the differential operator $\big(P_{\alpha  ,\beta
}(D)\big)^\gamma  $ is defined in the usual way, and we describe it
below for other $\gamma  $ in a manner similar to the way
$P(D)^\gamma  $ was described in (\ref{Eq4.9}) (for the special case
$\alpha  =\beta  =0).$  First, we recall the normalized Jacobi
polynomial $\varphi  _n$ given by
\begin{equation}\label{Eq10.5}
P_{\alpha  ,\beta  }(D)\varphi  _n = n(n+\alpha  +\beta  +1)\varphi  _n, \q
\int^1_{-1} \varphi  _n(x)\varphi  _k(x)w_{\alpha  ,\beta  }(x)dx =
\begin{cases}
1, &n=k\\
0, &n\ne k\end{cases}
\end{equation}
and the expansion of $f$ given by
\begin{equation}\label{Eq10.6}
f(x)\sim \sum^\infty  _{k=0} a_k \varphi  _k \q\text{\rm where} \q a_k
=\int^1_{-1} \varphi  _k(y)f(y)w_{\alpha  ,\beta  }(y)dy, \q P_k
f\equiv P^{(\alpha  ,\beta  )}_k f\equiv
a_k\varphi  _k\,.
\end{equation}

We now define $P_{\alpha  ,\beta  }(D)^\gamma  $ by
\begin{equation}\label{Eq10.7}
P_{\alpha  ,\beta  }(D)^\gamma   f\sim \sum^\infty  _{k=1}
\big(k(k+\alpha  +\beta  +1)\big)^\gamma  P_k f,
\end{equation}
with $P_k f$ given in (\ref{Eq10.6})
and $P_{\alpha  ,\beta  }(D)^\gamma  f \in B$ whenever there exists
$g\in B$ which satisfies $P^{(\alpha  ,\beta  )}_k g\equiv P_k g =
 \big(k(k+\alpha  +\beta  +1)\big)^\gamma  P_k f.$

In \cite[Theorem 7.1]{Da-Di05} it was shown for $1<p<\infty  $ that
\begin{equation}\label{Eq10.8}
\Vert  \varphi  ^r g^{(r)}\Vert  _{w_{\alpha  ,\beta  },p} \approx
\big\Vert  P_{(\alpha  ,\beta  )}(D)^{r/2} \big(g-S^{(\alpha  ,\beta
)}_{r-1}g\big)\big\Vert  _{w_{\alpha  ,\beta  },p}
\end{equation}
where $S^{(\alpha  ,\beta  )}_{r-1} g =\overset{r-1}{\us{k=0}{\sum}}
P^{(\alpha  ,\beta  )}_k g.$

Clearly, for $1<p<\infty  $ (\ref{Eq10.8}) implies
\begin{equation}\label{Eq10.9}
\begin{aligned}
\big[\text{\rm inf}\, \big(\Vert  f-g\Vert  _{w_{\alpha  ,\beta
},p} &+ t^r\Vert  \varphi  ^r g^{(r)}\Vert  _{w_{\alpha  ,\beta
},p}\big)\big] + t^r E_1(f)  _{w_{\alpha  ,\beta  },p}\\
& \approx\;\text{\rm inf}\, \big(\Vert  f-g\Vert  _{w_{\alpha  ,\beta
},p} + t^r \Vert  P_{(\alpha  ,\beta  )} (D)^{r/2}g\Vert
_{w_{\alpha  ,\beta  },p}\big).
\end{aligned}
\end{equation}
For $p=1$ and $p=\infty  $ (\ref{Eq10.9}) is not valid in the case
$\alpha  =\beta  =0.$

For smoothness given by $K  _\gamma  \big(f,P_{\alpha  ,\beta
}(D)^\gamma  ,t^{2\gamma  }\big)_{L_p(w)}$ many results related to
approximation were proved.  (Some are valid for all Banach spaces of
functions for which the Ces\`aro summability of some order of the
Jacobi expansion is bounded.)  The boundedness of the Ces\`aro summability of order
$r,$ $r>\max\, (\alpha  +\frac 12,\beta  +\frac 12)$ for
$L_{p,w}[-1,1]$ (with $w=w_{\alpha  ,\beta  })$
was given in \cite[p.~190]{Ch-Di97} as a corollary of
earlier results (see also \cite{Du-Xu} for the multivariate case),
that is
\begin{equation}\label{Eq10.10}
\Vert  C  ^r_nf\Vert  _{w_{\alpha  ,\beta  },p} \le C\Vert
f\Vert  _{w_{\alpha  ,\beta  },p},
\end{equation}
(where $C^r_nf \equiv C^r_n(f,x)$ is given by (\ref{Eq4.5}) with
$\varphi_k  $ of (\ref{Eq10.5})).  Therefore,
many theorems on polynomial approximation are valid.

The inequality (\ref{Eq10.10}) for $\alpha  >-1,$ $\beta  >-1$ and
$1\le p\le \infty  $ yields the following results:

\vs\noi
(A) A Bernstein-type inequality  given by
\begin{equation}\label{Eq10.11}
\Vert  P_{\alpha  ,\beta  }(D)^\gamma  P_n\Vert  _{w_{\alpha  ,\beta
},p} \le Cn^{2\gamma  }\Vert  P_n\Vert  _{w_{\alpha  ,\beta  },p}
\end{equation}
where $\gamma  >0$ (see \cite[(1.9)]{Ch-Di97} and \cite[(3.5)]{Di98}).

\vs
\noi
(B)  A direct or Jackson-type result  given by
\begin{equation}\label{Eq10.12}
E_n(f)_{w_{\alpha  ,\beta  },p} \le CK_\gamma  \big(f,P_{\alpha
,\beta  }(D)^\gamma  ,{n^{-2\gamma }}\big)_{w_{\alpha  ,\beta
},p}
\end{equation}
where $\gamma  >0$ (see \cite[(4.2)]{Ch-Di97} and \cite[(5.8) and
(5.22)]{Di98}).

\vs\noi
(C) A realization result  given by
\begin{equation}\label{Eq10.13}
K_\gamma  \big(f,P_{\alpha  ,\beta  }(D)^\gamma  ,
{n^{-2\gamma  }}\big)_{w_{\alpha  ,\beta  },p}
\approx \Vert  f-P_n\Vert  _{w_{\alpha  ,\beta  },p} +
n^{-2\gamma  }\Vert  P_{\alpha  ,\beta  }(D)^\gamma  P_n\Vert
_{w_{\alpha  ,\beta  },p}
\end{equation}
where $P_n$ is the best approximant to $f$ (i.e. satisfies $E_n (f)_{w_{\alpha
,\beta  },p} = \Vert  f-P_n\Vert  _{w_{\alpha  ,\beta  },p})$ or
$P_n = V_n f$ where $V_n$ is a de~la~Vall\'ee Poussin-type operator
(see \cite[(7.2)]{Di98}).  We remark that there are many operators
of the de~la~Vall\'ee Poussin type and we may choose
$$
V_n f\sim \sum^\infty  _{k=0}  \eta  \Big(\frac kn\Big) P_k
f\q\text{for} \q f\sim \sum^\infty  _{k=0}  P_k f
$$
where $\eta  (t) \in C^\infty  [0,\infty  ),$ $\eta  (t) =1$ for
$t\le 1,$ and $\eta  (t) =0$ for $t\ge 2$ for example.

\vs\noi
(D) The Marchaud-type inequality  given by
\begin{equation}\label{Eq10.14}
K_\gamma  \big(f,P_{\alpha  ,\beta  }(D)^\gamma  ,t^{2\gamma
}\big)_{w_{\alpha  ,\beta  },p}
\le Ct^{2\gamma  }\,\int^1_t\;\frac{K_\eta  \big(f,P_{\alpha  ,\beta
 }(D)^\eta  , u^{2\eta  })_{w_{\alpha  ,\beta  },p}}{u^{2\gamma  +1}}\; du, \q \eta  >\gamma
>0
\end{equation}
(see \cite[(5.25)]{Ch-Di97} and \cite[(6.7)]{Di98}).

\vs\noi
(E)  The converse-type inequality  given by
\begin{equation}\label{Eq10.15}
K_\gamma  \Big(f,P_{\alpha  ,\beta  }(D)^\gamma  ,
\frac{1}{n^{2\gamma  }}\Big)_{w_{\alpha  ,\beta  },p}
\le Cn^{-2\gamma  } \sum^n_{k=1} k^{2\gamma  -1}E_k(f)_{w_{\alpha
,\beta  },p}
\end{equation}
(see \cite[(5.23)]{Ch-Di97} and \cite[(6.6)]{Di98}).

\vs\noi
(F) As a result of simultaneous approximation, one has
\begin{equation}\label{Eq10.16}
E_n (f)_{w_{\alpha  ,\beta  },p} \le Cn^{-2\gamma
}E_n\big(P_{\alpha  ,\beta  }(D)^\gamma  f\big)_{w_{\alpha  ,\beta
},p}
\end{equation}
whenever $P_{\alpha  ,\beta  }(D)^\gamma   f\in L_{w_{\alpha  ,\beta
 },p} $ (see \cite[(7.3)]{Di98}).

In addition, we have the sharp Marchaud-type inequality
\begin{equation}\label{Eq10.17}
K_\gamma  \big(f,P_{\alpha  ,\beta  }(D)^\gamma  ,t^{2\gamma
}\big)_{w_{\alpha  ,\beta  },p}
\le Ct^{2\gamma  }\,\Big\{\int^1_t\;\frac{K_\eta  \big(f,P_{\alpha
,\beta  }(D)^\eta  ,u^{2\eta  })^q_{w_{\alpha  ,\beta  },p}} {u^{2\gamma  q+1}}
\,du \Big\}^{1/q}
\end{equation}
for $0<\gamma  <\eta  ,$ $1<p<\infty  $ and $q=\min\, (p,2),$ which
was proved in \cite[(6.8)]{Da-Di05} and we also have the corresponding
sharp converse result
\begin{equation}\label{Eq10.18}
K_\gamma  \big(f,P_{\alpha  ,\beta  }(D)^\gamma  ,t^{2\gamma
}\big)_{w_{\alpha  ,\beta  },p}
\le Ct^{2\gamma  }\Big(\sum_{1\le n<1/t}\, n^{2\gamma  q-1}
E_n(f)^q_{w_{\alpha  ,\beta  },p} \Big)^{1/q}
\end{equation}
for $0<\gamma  ,$ $1<p<\infty  $ and $q=\,\min\,(p,2)$ (see
\cite[(6.10)]{Da-Di05}).

Another inequality related to Jacobi weights is the Nikol'skii-type
inequality
\begin{equation}\label{Eq10.19}
\Vert  P_n\Vert  _{w_{\alpha  ,\beta  },q} \le Cn^{\gamma  (\frac
1p-\frac 1q)} \Vert  P_n\Vert  _{w_{\alpha  ,\beta  },p}
\q\text{for} \q P_n\in \Pi_n,
\end{equation}
where $\gamma  =\,\max\,\big(2+2\,\max\,(\alpha  ,\beta  ),1\big)$
and $0<p<q\le \infty  $ (see \cite[Theorem 6.6]{Di-Ti05} for a
somewhat more general result and a simple proof).  The inequality
(\ref{Eq10.19}) is also best possible (like (\ref{Eq9.14}) and
(\ref{Eq9.15})) as equality holds for $p=2$ and $q=\infty  .$

Recently (see \cite[Theorem 6.1]{Da-Di-Ti}) a sharp version of the
Jackson inequality (\ref{Eq10.12}) for $1<p<\infty  $ $(p\ne
1,\infty  )$ was given by
\begin{equation}\label{Eq10.20a}
2^{-2n\gamma  } \Big\{\sum^n_{j=1} 2^{2j\gamma  s}
E_{2^j}(f)^s_{L_{p,w_{\alpha  ,\beta  }}[-1,1]} \Big\}^{1/s}
\le CK_\gamma  \big(f,P_{\alpha  ,\beta  }(D)^\gamma  ,2^{-2n\gamma
}\big)_{L_{p,w_{\alpha  ,\beta  }}[-1,1]}
\end{equation}
where $s=\max\,(p,2)$ and $\gamma  >0.$  Similarly, a form
equivalent to (\ref{Eq10.20a}) comparing $K$-functionals and
extending (\ref{Eq10.14}) for $1<p<\infty  $ was achieved (see
\cite[(6.3)]{Da-Di-Ti}) and is given for $\zeta  >\gamma  $ and
$s=\max\,(p,2)$ by
\begin{equation}\label{Eq10.21a}
\begin{aligned}
2^{-nr} \Big\{\sum^n_{j=1} &2^{2j\gamma  s} K_\zeta  \big(f,P_{\alpha
 ,\beta  } (D)^\zeta  ,2^{-2j\zeta  }\big)^s_{L_{p,w_{\alpha  ,\beta
 }}[-1,1]} \Big\}^{1/s}\\
&\le CK_{\gamma}  \big(f,P_{\alpha  ,\beta  }(D)^\gamma  ,2^{-2n\gamma
}\big)_{L_{p,w_{\alpha  ,\beta  }}[-1,1]}.
\end{aligned}
\end{equation}

For Jacobi weights the Durrmeyer operator
\begin{equation}\label{Eq10.20}
M^{(\alpha  ,\beta  )}_n f\equiv M^{(\alpha  ,\beta  )}_n(f,x)
=\sum^n_{k=1} \big(A^{(\alpha  ,\beta
)}_{n,k}\big)^{-1}P_{n,k}(x)\int^1_{0}P_{n,k}(y)f(y)w_{\alpha
,\beta  }(y)dy
\end{equation}
where $P_{n,k}$ is given in (\ref{Eq8.1}) and $ A_{n,k}
=\int^{1}_{0} P_{n,k}(y)w_{\alpha  ,\beta  }(y)dy$
satisfies a strong converse inequality with the $K\text{\rm
-functional}$ given in (\ref{Eq10.3}).
That is,
\begin{equation}\label{Eq10.21}
\Vert  f-M^{(\alpha  ,\beta  )}_n f\Vert  _p \approx
K_1\Big(f,P_{\alpha  ,\beta  }(D),\,\frac 1n\Big)_p
\end{equation}
where $K_1\big(f,P_{\alpha  ,\beta  }(D),\,\frac 1n\big)_p$ is given
in (\ref{Eq10.3}) (with $\gamma  =1$ and $t^2=n^{-1}).$
This will be discussed, together with
its multivariate analogues, in Section \ref{Sec12}.

For more information on weighted approximation with Jacobi weights
see Section \ref{Sec18}.

\section{Weighted moduli, Freud weights}\label{Sec11}

To approximate functions on ${\IR}$ by polynomials, one needs to consider
weighted approximation.  A detailed discussion of this problem
appears in the major survey on that topic by Lubinsky \cite{Lu}.
Other important sources on the subject are the books by Levin and
Lubinsky \cite{Le-Lu} and by Mhaskar \cite{Mh}.  Here we just
briefly outline the results related to polynomial approximation and
put them in the context of this survey.

To investigate the rate of approximation by polynomials to a
function in $L_{p}(W,\IR)$ given by the norm or quasinorm $\Vert
Wf\Vert  _{L_p(\IR)},$
one must first ascertain for which type of
weights $W(x)$  polynomials are dense in $L_{p}(W,\IR).$  Necessary
and sufficient conditions on $W(x)$ were given by Akhieser,
Carleson, Mergelian and Pollard.  (For a more detailed discussion
see \cite[Section 1]{Lu}).

We deal here with Freud weights (see \cite{Fr}, \cite{Lu} and
\cite{Mh}) which are given by
$W(x) =$
\newline $\exp\,\big(-Q(x)\big)$ with $Q(x)$ an even continuous
function, with
$Q^{\pr\pr}(x)$ continuous, $Q^\pr(x)$ positive in $(0,\infty  ),$ and
for some $a,b>0$
\begin{equation}\label{Eq11.1}
a\le \frac{xQ^{\pr\pr}(x)}{Q^\pr(x)} \le b\q\text{\rm for}\q x\in
(0,\infty  ).
\end{equation}
In fact, the results are valid for somewhat more general $Q(x),$ but
the most prominent cases, that is, when $Q(x) = \vert  x\vert
^\alpha  ,$ $\alpha  >1$ already satisfy the above conditions. To
define moduli of smoothness, $K$-functionals, and realization
functionals for the spaces $L_{p,W}(\IR)$ of functions satisfying
$Wf\in L_p(\IR),$
one needs to define the
Mhaskar-Rahmanov-Saff number $a_n$ which is the root of
\begin{equation} \label{Eq11.2}
n =\frac 2\pi \;\int^1_0 \;\frac{a_n tQ^\pr(a_nt)dt}{\sqrt{1-t^2}}
\,, \q n>0.
\end{equation}
The number $a_n$ gives rise to the crucial Remez-type inequalities
\begin{equation}\label{Eq11.3}
\Vert  PW\Vert  _{L_\infty  (\IR)} \le \Vert  PW\Vert  _{L_\infty
(-a_n,a_n)} \q\text{\rm for}\q  P\in \Pi_n,
\end{equation}
and for $0<p<\infty  $
\begin{equation}\label{Eq11.4}
\Vert  PW\Vert  _{L_p(\IR)} \le (1+e^{-Cn}) \Vert  PW\Vert
_{L_p(\vert  x\vert  <a_n+\varepsilon  )} \q\text{\rm for}\q  P\in \Pi_n
\end{equation}
where $C=C(p,\varepsilon  ,W)$ does not depend on $n.$

The modulus of smoothness is given by
\begin{equation}\label{Eq11.5}
\omega  ^r(f,W,t)_p =\us{0<h\le t}\sup\; \Vert  W\Delta  ^r_h f\Vert
 _{L_p[-\sigma  (h),\sigma  (h)]} + \us{P\in \Pi_r}\inf\;\Vert
W(f-P)\Vert  _{L_p[\vert  x\vert  \ge \sigma  (t)]}
\end{equation}
where
\begin{equation}\label{Eq11.6}
\sigma  (h) =\,\inf\,\Big\{a_n:\,\frac{a_n}{n}\le h\Big\}
\end{equation}
and
$$
\Delta  ^r_h f(x) = \sum^r_{i=0} \,\binom r i (-1)^i
f\big(x+\frac{rh}{2} - ih\big).
$$
We observe that for $W_\alpha  (x) = e^{-\vert  x\vert  ^\alpha  }$
with
$\alpha  >1,$ $a_n\approx n^{1/\alpha  }$ and $\sigma
\big(\frac{a_n}{n}\big)\approx a_n.$  We also note that, following
\cite{Di-Lu}, the second term on the right of (\ref{Eq11.5}) is
different from that in \cite[Chapter 11, (11.2.2)]{Di-To87} to
accommodate the space $L_{p,W}(\IR)$
of functions for which $Wf\in L_p(\IR)$
with $0<p<1.$  (For $1\le p\le
\infty  $ the two forms are equivalent.)

The rate of best polynomial approximation is given by
\begin{equation}\label{Eq11.7}
E_n(f)_{W,p} =\,\inf\,(\Vert  W(f-P)\Vert  _{L_p(R)}: P\in \Pi_n).
\end{equation}

The $K$-functional is given by
\begin{equation}\label{Eq11.8}
K_r(f,W,t^r)_p =\,\us g\inf\, \big(\Vert  (f-g)W\Vert  _{L_p(\IR)} +
t^r\Vert  g^{(r)}W\Vert  _{L_p(\IR)}\big),
\end{equation}
which, following the technique in \cite{Di-Hr-Iv}, satisfies
\begin{equation}\label{Eq11.9}
K_r(f,W,t^r)_p =0 \q \text{\rm for} \q Wf\in L_p(\IR) \q \text{\rm
when}\q  0<p<1.
\end{equation}
The realization functional $\wt{R}_r(f,W,t^r)_p$ is given by
\begin{equation}\label{Eq11.10}
\wt{R}_r (f,W,t^r)_p = \inf\Big\{\Vert  (f-P)W\Vert  _{L_p(\IR)}
+ t^r\Vert  P^{(r)}W\Vert  _{L_p(\IR)}: P\in \Pi_n,
 \; n = \inf\,\big(k: \frac{a_k}{k}\le t\big)\Big\}
\end{equation}
(see \cite{Di-Lu}), or by its equivalent form
\begin{equation}\label{Eq11.11}
R_r(f,W,t^r)_p = \Vert  (f-P_n)W\Vert  _{L_p(\IR)} + t^r \Vert
P^{(r)}_nW\Vert  _{L_p(\IR)}
\end{equation}
where
$ n = \,\inf\,\big(k:\,\frac{a_k}{k}\le t\big)$ for $P_n$
satisfying $E_n(f)_{W,p} = \Vert  (f-P_n)W\Vert  _p,$ and $P_n\in \Pi_n.$

For $1\le p\le \infty  $ a de la Vall\'ee Poussin-type expression
$V_nf$ can replace $P_n$ in (\ref{Eq11.11}).  We note that because
of the boundedness of the Ces\`aro summability of order $1,$ $V_n$
may be given by the classical form
\begin{equation}\label{Eq11.12}
V_n(f,x) = 2C_{2n}(f,x) - C_n(f,x), \q C_n(f,x) = C^1_n(f,x)
\end{equation}
where $C^\ell_n(f,x)$ is given by (\ref{Eq4.5})
$\big(C_n(f,x) = C^1_n(f,x)\big)$
with $\varphi
_k\in \Pi_{k+1}\,,$  orthonormal polynomials with respect to $W,$ i.e.
$$
\int_{\IR} \varphi  _k\varphi  _\ell W^2=\begin{cases} 1, &k= \ell\\ 0, &k\ne
\ell.\end{cases}
$$
The properties $V_nf\in \Pi_{2n},$ $V_n P = P$ for $P\in \Pi_n$ and
$\Vert W V_nf\Vert  _p \le A\Vert W f\Vert  _p$ are clear and given
for instance in \cite[p.~70]{Mh}.  We note that $\eta  _n f$ given
by (\ref{Eq5.4}) is also a de~la~Vall\'ee Poussin-type operator, and
$\Vert  W\eta  _n f\Vert  _p\le C\Vert W f\Vert  _p$ follows from
$\Vert  WC_n f\Vert  _p\le A_1\Vert W f\Vert  _p$ and summation by
parts (Abel transformation).

We define the realization with $V_n$ by
\begin{equation}\label{Eq11.13}
R^*_r(f,W,t^r)_p = \Vert  (f-V_nf)W\Vert  _{L_p(\IR)} + t^r\,\Big\Vert
W\Big(\frac{d}{dx}\Big)^r V_nf\Big\Vert  _{L_p(\IR)}\,,
\end{equation}
where $n= \inf\,(k:\frac{a_k}{k}\le t)$
and we note that as $V_nf$ and $R^*_r$ are not defined for
$0<p<1,$ the equivalence
$R^*_r(f,W,t^r)_p \approx R_r(f,W,t^r)_p$ is valid only for
$1\le p\le \infty  .$

For $0<p\le \infty  $ one has
\begin{equation}\label{Eq11.13a}
\omega  ^r(f,W,t)_p\approx R_r(f,W,t^r)_p\approx
\wt{R}_r(f,W,t^r)_p\,.
\end{equation}
For $1\le p\le \infty  $ one has
\begin{equation}\label{Eq11.14}
\omega  ^r(f,W,t)_p\approx K_r(f,W,t^r)_p \approx R^*_r(f,W,t^r)_p\,.
\end{equation}

The Markov-Bernstein inequality is given by
\begin{equation}\label{Eq11.15}
\Vert  P^\pr W\Vert  _{L_p(\IR)} \le C\,\frac{n}{a_n}\, \Vert  PW\Vert
 _{L_p(\IR)}, \q P\in \Pi_n, \q 0<p\le \infty
\end{equation}
where $C=C(p,W)$ does not depend on $n$ or $P$ (see the discussion in
\cite{Le-Lu}).  The Jackson inequality (see \cite[p.~102]{Di-Lu}) is
given by
\begin{equation}\label{Eq11.16}
E_n (f) _p\le C_1 \omega  ^r\big(f,W,C_2(a_n/n)\big)_p, \q 0<p\le \infty  ,
\end{equation}
where $C_i$ are independent of $n$ and $f.$ Also $C_2$ can be
replaced by $1$ for $1\le p\le \infty  $ (see \cite[p.~104]{Di-Lu}).
 Furthermore, we have the converse (to (\ref{Eq11.16})) result (see
\cite[p.~105]{Di-Lu}) given by
\begin{equation}\label{Eq11.17}
\omega  ^r(f,W,t)^q_p \le C\Big(\frac{a_n}{n}\Big)^{rq}
\,\sum^\ell_{j=0}\,\Big(\frac{2^j}{a_{2^j}}\Big)^{rq}
E_{2^j}(f)^q_{w,p}, \q 0<p\le \infty  ,\q q=\,\min\,(p,1)
\end{equation}
for $t$ small enough and $n=\,\inf\,\big(k:\frac{a_k}{k}\le t\big).$
As a corollary of (\ref{Eq11.16}) and (\ref{Eq11.17}), we have the
Marchaud-type inequality (see \cite[p.~105]{Di-Lu}) for $0<p\le \infty
,$ $q=\,\min\,(p,1)$ given by
\begin{equation}\label{Eq11.18}
\omega  ^r(f,W,t)_p \le C_1 t^r\,\Big\{\int^1_t\;\frac{\omega
^{r+1}(f,W,t)^q_p}{u^{rq+1}}\;du +\Vert  fW\Vert  ^q_p\Big\}^{1/q}.
\end{equation}
For $1\le p\le \infty  $ one has
\begin{equation}\label{Eq11.19}
\omega  ^r(f,W,t)_p\le Ct^r\Vert  f^{(r)}W\Vert  _{L_p(\IR)}\,.
\end{equation}
As the saturation class of $\omega  ^r(f,W,t)_p$ for $0<p<1$ is
$O(t^\gamma  )$ with $\gamma  >r,$ (\ref{Eq11.19}) is not useful for
that range.

I conjecture that for $1<p<\infty  $ a sharp Marchaud and a sharp
Jackson inequality will eventually be established.  That is, for
$1<p<\infty , $ (\ref{Eq11.17}) and (\ref{Eq11.18}) will be proved
with $q=\,\min\,(p,2)$ rather than with $q=\,\min\,(p,1),$ and
an analogue of (\ref{Eq2.8}) with
$s=\,\max\,(p,2)$ will replace (\ref{Eq11.16}).

I will deal with Nikol'skii and Ul'yanov-type inequalities in Section
\ref{Sec13} and with multivariate analogues (essentially the lack
thereof) in Section \ref{Sec12}.

\section{Multivariate polynomial
approximation}\label{Sec12}

The space of polynomials of total degree smaller than $n,\Pi_n$ is
given by
\begin{equation}\label{Eq12.1}
\Pi_n =\;\spa\,\Big\{x^{\alpha  _1}_1 \cdots x^{\alpha  _d}_d: \,
\alpha  _i=0,1,2,\dots, \,\alpha _1+\dots +\alpha  _d <n\Big\}.
\end{equation}

It is a natural question to ask for what spaces of functions and on
what domains one can extend the Bernstein, Jackson, Marchaud and
other inequalities.  In this section we will outline the progress
made after the text \cite{Di-To87} appeared.

For a convex set $S$ in $\IR^d$ it was shown in \cite[Theorem2.1]{Di92}
that the Bernstein inequality on the interval can be copied to read
for $0<p\le \infty  $ and $r=1,2,\dots$
\begin{equation}\label{Eq12.2}
\Big\Vert  \wt d(\bx,\bxi  )^{r/2}\Big(\pd{}{\bxi  }\Big)^r
P_n(\bx)\Big\Vert  _{L_p(S)} \le Cn^r \Vert  P_n\Vert  _{L_p(S)}
\q\text{\rm for}\q P_n\in \Pi_n
\end{equation}
with $C$ independent of $\bxi  ,n,P_n$ and $S$ and with $\wt d(\bx,\bxi
)$ which was introduced in \cite{Di92} and given by
\begin{equation}\label{Eq12.3}
\varphi  _{\bxi}  (\bx)^2\equiv \wt d(\bx,\bxi  ) \equiv \wt d_S(\bx,\bxi  ) =
\us{\bx+\lambda  \bxi  \in S}{\sup\,\lambda} \;  \us{\bx-\mu  \bxi  \in
S}{\sup\, \mu}  , \q \bx\in S, \q \vert  \bxi  \vert  =1.
\end{equation}

We note that when $S$ is unbounded, (\ref{Eq12.2}) is meaningless,
and the same is true when the interior of $S$ is empty and $p<\infty
 ,$ so we may as well consider only bounded convex sets $S$ with
non-empty interior. The introduction of $\varphi  _{\bxi}  (\bx)^2 =
\wt d(\bx,\bxi)$ in \cite{Di92}, which in (\ref{Eq12.2}) yields a
constant independent of $S,$ is natural since for $S=[-1,1]\subset
\IR,$ (for which only $\xi  $ is equal to $\pm e$ where $e= (0,1)$
is possible), $\varphi  (x) =\wt d_S(x,\xi  )^{1/2} = \sqrt{1-x^2}\,.$

The Markov-type inequality for a bounded convex set with non-empty
interior $S\subset R^d$ was given by (see \cite[Theorem 4.1]{Di92})
\begin{equation}\label{Eq12.4}
\Big\Vert  \Big(\pd{}{\bxi  }\Big)^r P_n \Big\Vert  _{L_p(S)} \le
Cn^{2r}\Vert  P_n\Vert  _{L_p(S)}, \q P_n \in \Pi_n
\end{equation}
where $C$ depends on $S$ and $\bxi$ but not on $n$ or $P_n.$

The Remez-type inequality (see \cite[Theorem 3.1]{Di92}) for a
bounded convex set $S$ with non-empty interior is given by
\begin{equation}\label{Eq12.5}
\Vert  P_n\Vert  _{L_p(S)} \le C(p,L,S)\Vert  P_n\Vert
_{L_p(S(L,n))},
\end{equation}
where $S(L,n) =\{\bu:B(\bu,L/n^2)\subset S\}$ and $B(\bx,r)$ is the
ball of center $\bx$ and radius $r.$

Inequalities like (\ref{Eq12.4}) were studied extensively and for
various more general multivariate domains, but the polynomial
approximation and its relations to concepts of smoothness
generalizing $\omega
 ^r_\varphi  (f,t)$ were not tackled even if one assumes that we are
dealing with a general bounded convex set with non-empty interior.

In \cite[Chapter 12]{Di-To87} the direct and the weak converse
inequalities were proved for $L_p(S)$ (when $1\le p\le \infty  )$
and where $S$ is a simple polytope.  We recall that a polytope is
the convex hull of finitely many points in $\IR^d,$ and a simple
polytope is a polytope all of whose vertices are joined to other
vertices by exactly $d$ edges.  A simplex and a box or a cube are
perhaps the most familiar simple polytopes. The Egyptian pyramid is
not a simple polytope.

The moduli of smoothness on a polytope $S$
can be given by
\begin{equation}\label{Eq12.6}
\wt\omega  ^r_S(f,t)_{L_p(S)} =
\us{ {{\vert  h\vert  <t}\atop{\bu \in E(S)}}\atop{\bv = \bu/\vert
\bu\vert   }}{\text{\rm sup}}
\,\Vert  \Delta  ^r_{h
\wt d(\bx,\bv)^{1/2}\bv}f(\bx)\Vert  _{L_p(S)}
\end{equation}
where $E(S)$ is the set of edges of $S$ and
\begin{equation}\label{Eq12.7}
\Delta  ^r_{h\wt d(\bx,\bv)^{1/2}\bv} f(\bx)
= \begin{cases} \os r{\us{k=0}\sum} (-1)^k\binom rk f\Big(\bx + \big(\frac
r2-k\big)h\wt d(\bx,\bv)^{1/2}\bv\Big), &\bx\pm\frac r2 h\wt
d(\bx,\bv)^{1/2}\bv\in S\\
0\;, &\text{\rm otherwise}.\end{cases}
\end{equation}

We may also define $\omega  ^r_S(f,t)_p$ as
\begin{equation}\label{Eq12.8}
\omega  ^r_S(f,t)_{L_p(S)} =\us{ {\vert  h\vert  \le
t}\atop{\vert  \bv\vert  =1 \; \bv \in \IR^d}}\sup \;\Vert
\Delta  ^r_{h\wt d(\bx,\bv)^{1/2}\bv} f(\bx)\Vert  _{L_p(S)}\,,
\end{equation}
which is defined for all convex sets $S.$  It is known that for
$p=1$ and $p=\infty  $ and for a simple polytope $S,$  $\wt \omega
_S(f,t)_p$ is not equivalent to $\omega  _S(f,t)_p$ (see
\cite[Remark 12.2.1]{Di-To87}).  I expect that for $1<p<\infty  $ and a
simple polytope $S,$ $\omega  _S(f,t)_p\approx \wt\omega  _S(f,t)_p,$
though this was not yet proved.  The definition of the moduli of
smoothness in (\ref{Eq12.6}) and (\ref{Eq12.8}) are different from
those in \cite{Di-To87} only in style, using here the somewhat more
convenient concept $\wt d(\bx,\bxi  )$ given in (\ref{Eq12.3}).  For a
simple polytope $S$ the Jackson inequality, given by
\begin{equation}\label{Eq12.9}
E_n(f)_{L_p(S)} \equiv \us{P\in \Pi_n}{\text{\rm inf}}\;\Vert
f-P\Vert  _{L_p(S)} \le C\wt \omega  ^r_S\Big(f,\frac
1n\Big)_{L_p(S)}, \q 0<p\le \infty
\end{equation}
was proved in \cite[Chapter 12]{Di-To87} for $1\le p\le \infty  $ and
in \cite[Theorem 1.1]{Di96} for $0<p<1.$

If the boundary effect is ignored, a Jackson-type estimate is
possible for a much more general domain $\SD.$  We define the
modulus $\omega  ^r(f,t)_{L_p(D)}$ by
\begin{equation}\label{Eq12.10}
\omega  ^r(f,t)_{L_p(\cal D)} =\us{\vert  \bh\vert  \le t}\sup\, \Vert
\Delta  ^r_{\bh}f\Vert  _{L_p(\cal D)}
\end{equation}
where
\begin{equation}\label{Eq12.11}
\Delta  ^r_{\bh} f(\bx) =\begin{cases}
\os r{\us{k=0}\sum} (-1)^k \binom rk f\Big(\bx +\bh \big(\frac r2
-k\big)\Big), &\Big\{\bx + \tau  \bh:\vert  \tau  \vert  \le \frac
r2\Big\}\subset \cal D\\
0, &\text{\rm otherwise.}
\end{cases}
\end{equation}

If one can extend $f$ to be defined on a cube $Q$ in such a way that
$g= Ef$ defined on $Q$ satisfies
\begin{equation}\label{Eq12.12}
\omega  ^r(g,t)_{L_p(Q)} \le C\omega  ^r(f,t)_{L_p(D)} \q\text{\rm
and}\q
f(\bx)
=g(\bx)\q\text{\rm for}\q \bx\in \cal D,
\end{equation}
and if $Q\supset Q^*\supset D$ where $Q^*+\bxi \subset Q$ for $\vert
 \bxi \vert  \le 1,$ then (\ref{Eq12.9}) clearly implies
\begin{equation}\label{Eq12.13}
E_n (f) _{L_p(\cal D)} \le C\omega  ^r(f,t)_{L_p(\cal D)}.
\end{equation}
Such an extension of $f$ on $\cal D$ is discussed for instance in
\cite{De-Sh}, and the results there are given for all $p$ and for
many bounded domains.

However, formulae like (\ref{Eq12.13}) are a departure from the topic
of this survey, as the effect of being near the boundary is
neglected.  Moreover, there is no hope of having a matching converse
result to (\ref{Eq12.13}) as the results described below will imply.

For a simple polytope $S$ the converse result (see \cite[Theorem
12.2.3]{Di-To87} for $1\le p\le \infty  $ and \cite{Di96} for
$0<p<1)$ is given by
\begin{equation}\label{Eq12.14}
\begin{aligned}
\wt \omega  ^r_S(f,t)_{L_p(S)} &\le \omega  ^r_S(f,t)_{L_p(S)}\\
&\le Mt^r \Big\{\us{1\le k\le \frac 1t}\sum\,k^{rq-1}
E_k(f)^q_{L_p(S)}\Big\}^{1/q}
\end{aligned}
\end{equation}
with $q=\min\,(p,1).$

The Marchaud-type variation of (\ref{Eq12.14}) given by
\begin{equation}\label{Eq12.15}
\omega  ^r_S(f,t)_{L_p(S)} \le
C\Big\{\int^1_t\,\frac{\omega  ^{r+1}_S(f,u  )^q_p }{u^{rq+1}}\;du +
\Vert  f\Vert  ^q_{L_p(S)}\Big\}^{1/q}
\end{equation}
for $0<p\le \infty  ,$ $q=\min \,(p,1)$ and $C=C(r,p,S)$ independent
of $f$ and $t$ was proved for a simple polytope $S$ and $1\le p\le
\infty  $ in \cite{Di-To87} and for $0<p<1$ in \cite[Theorem
5.1]{Di96}.

We note that one does not have the sharp versions of the Jackson,
Marchaud and the (weak) converse inequality, that is, the
analogues of (\ref{Eq2.8}), (\ref{Eq6.6}) and (\ref{Eq6.4}) to
replace (for simple polytopes and $1<p<\infty  )$ the inequalities
(\ref{Eq12.9}), (\ref{Eq12.14}) and (\ref{Eq12.15}) respectively. I
believe that such results will eventually be proved.  For the
simplex (see \cite{Da-Di07}) (\ref{Eq12.10}) was extended, and that
result will be described later in this section.

For $1\le p\le \infty  $ and a simple polytope $S,$ $\wt \omega
^r_S(f,t)_{L_p(S)}$ and $\omega  ^r_S(f,t)_{L_p(S)} $ are equivalent
to the $K$-functionals $\wt K_{r,S}(f,t^r)_{L_p(S)}$ and
$K_{r,S}(f,t^r)_{L_p(S)}$ respectively as stated in the following
formulae:
\begin{equation}\label{Eq12.16}
\begin{aligned}
\wt K_{r,S}(f,t^r)_{L_p(S)}
&\equiv \,\us g{\text{\rm inf}}\;\Big(\Vert  f-g\Vert  _{L_p(S)} +
t^r \,\us{\bxi  \in E(S)}\sup\;\Big\Vert  \varphi  ^r_{\bxi }
\Big(\pd{}{\bxi  }\Big)^r g\Big\Vert  _{L_p(S)}\Big)\\
&\approx \wt \omega  ^r_S(f,t)_{L_p(S)}
\end{aligned}
\end{equation}
and
$$
\begin{aligned}
K_{r,S}(f,t^r)_{L_p(S)}
&\equiv \;\us g{\text{\rm inf}}\;\Big(\Vert  f-g\Vert  _{L_p(S)} +
t^r \,\us \bv{\text{\rm sup}}\, \Big\Vert  \varphi  ^r_{\bv}
\Big(\pd{}{\bv}\Big)^r g\Big\Vert  _{L_p(S)}\Big)\\
&\approx \omega  ^r_S(f,t)_{L_p(S)}
\end{aligned}\eqno{(12.16)^\prime}
$$
where $\varphi  _\xi  $ is given by (\ref{Eq12.3}).

While it was not shown explicitly in \cite{Di-Hr-Iv}, the method
there for $f\in L_p(S)$ with $0<p<1$ implies the equality
\begin{equation}\label{Eq12.17}
K_{r,S}(f,t^r)_{L_p(S)} = \wt K_{r,S}(f,t^r)_{L_p(S)} = 0.
\end{equation}
This just adds to the interest in the realization concept.

It was shown in \cite[(4.4)]{Di96} that
\begin{equation}\label{Eq12.18}
\wt \omega  ^r_S \big(f,\frac 1n\big)_{L_p(S)}
\approx \wt R_{r,S}(f,n^{-r})_{L_p(S)} \equiv \Vert  f-P_n\Vert  +
n^{-r}\,\us{\bxi  \in E_S}{\text{\rm sup}}\, \Big\Vert  \varphi
^r_{\bxi}  \Big(\pd{}{\bxi  }\Big)^r P_n\Big\Vert  _{L_p(S)}
\end{equation}
where $0<p\le \infty  ,$ $r=1,2,\dots,$ $S$ is a simple polytope and
$P_n$ a near best $n$-th degree polynomial approximant of $f$ in
$L_p(S).$  (A similar result holds for $\omega  ^r_S\big(f,\frac
1n\big)_{L_p(S)}.)$

The multivariate Bernstein polynomials on a simplex $S\subset \IR^d$ where
\begin{equation}\label{Eq12.19}
S=\Big\{(x_1,\dots,x_d):\, 0\le x_i,\, 0\le x_0 = 1 -\sum^d_{i=1}
x_i\Big\}
\end{equation}
is given by
\begin{equation}\label{Eq12.20}
B_n(f,\bx) = \sum_{\vert  \bk\vert  \le n} P_{n,\bk} (\bx)
f\Big(\frac \bk n\Big)
\end{equation}
where $\bk = (k_1,\dots,k_d), $\; $\vert  \bk\vert  =\os
d{\us{i=1}\sum} k_i,$\; $k_0= 1-\vert  \bk\vert  ,$ \; $x_0 = 1-\vert
 \bx\vert  ,$ and
\begin{equation}\label{Eq12.21}
P_{n,\bk} (\bx) =\frac{n!}{k_0! k_1! \dots k_d!} \, x^{k_0}_0
x^{k_1}_1 \dots x^{k_d}_d\,.
\end{equation}

The Bernstein polynomial is a contraction operator on $C(S)$ with
the Voronovskaja
\begin{equation}\label{Eq12.22}
n\big(B_n f(\bx) - f(\bx)\big) \to \frac 12\;\sum_{\bxi  \in
E_S,\,\bxi ={\pmb\zeta}  /\vert {\pmb \zeta}  \vert  } \wt
d_S(\bx,\bxi)\Big(\pd{}{\bxi  }\Big)^2 f(\bx) \equiv P_S(D).
\end{equation}
In earlier texts a long form of (\ref{Eq12.22}) was discussed, but
the use of (\ref{Eq12.3}) yields the compact expression
(\ref{Eq12.22}), which also demonstrates the intrinsic symmetry among
the edges.

The strong converse inequality
\begin{equation}\label{Eq12.23}
\Vert  B_n f-f\Vert  _{C(S)} \approx
K\big(f,P_S(D),n^{-1}\big)_{C(S)}
=\us g{\text{\rm inf}}\, \big(\Vert  f-g\Vert  _{C(S)} + n^{-1}\Vert
 P_S(D)g\Vert  _{C(S)}\big)
\end{equation}
for $S,$ $B_nf\equiv B_n(f,x)$ and $P_S(D)$ given in
(\ref{Eq12.19}), (\ref{Eq12.20}) and (\ref{Eq12.22}) respectively
was claimed in \cite{Zh}.

For $\alpha  <2$ one has
$$
E_n(f)_{C(S)} = O\Big(\frac{1}{n^\alpha  }\Big) \Longleftrightarrow \Vert
 B_nf-f\Vert  _{C(S)} = O\Big(\frac{1}{n^{\alpha  /2}}\Big)\,,
$$
which was known earlier.

The multivariate Durrmeyer-Bernstein polynomial approximation on the
simplex $S$ (given by (\ref{Eq12.19})) with the Jacobi weight
$w_{\bal}(\bx)$ given by
\begin{equation}\label{Eq12.24}
w_{\bal}(\bx) = \big(1-\vert  \bx\vert  \big)^{\alpha  _0} x^{\alpha
_1}_1 \dots x^{\alpha  _d}_d, \q \alpha  _i >-1, \q \bx \in S
\end{equation}
for $\bal = (\alpha  _0,\alpha  _1,\dots,\alpha  _d),$ $\bx =
(x_1,\dots,x_d  )$ and $x_0 = 1-\vert  \bx\vert  $ is defined by
\begin{equation}\label{Eq12.25}
M_{n,\bal}(f,\bx) = \sum_{\vert  \bk\vert  \le n} P_{n,\bk} (\bx)
A^{-1}_{n,\bk,\bal} \int_S P_{n,\bk}(\by) f(\by) w_{\bal}(\by)d\by
\end{equation}
where $\int_S P_{n,\bk}(\by) w_\alpha  (\by) dy = A_{n,\bk,\bal}\,.$

The behaviour of $M_{n,\bal}(f,\bx)$ and its rate of approximation
in $L_{p,W_{\bal}  }(S)$ were studied in many articles (see
\cite{Ch-Di-Iv}, \cite{De85}, \cite{Di95}, \cite{Zh} and others).

The Voronovskaja of $M_{n,\bal}f$ is given by
\begin{equation}\label{Eq12.26}
\begin{aligned}
n\big(M_{n,\bal}f(\bx) - f(\bx)\big) &\to \frac 12\, \sum_{\xi  \in
E_S} \,\frac{1}{w_{\bal}  (\bx)}\; \pd{}{\bxi  }\; \wt d_S(\bx,\bxi  )
w_{\bal}  (\bx)\;\pd{}{\bxi  }\; f(\bx)\\
&= \frac 12\; \wt P_{S,\bal  }(D),
\end{aligned}
\end{equation}
which clearly exhibits both its self-adjointness and dependence on
$w_{\bal}(\bx).$  In \cite{Di95} the technique of Knoop and Zhou is
used (and modified) to obtain the strong converse inequality
\begin{equation}\label{Eq12.27}
\begin{aligned}
\Vert  M_{n,\bal  }f-f\Vert  _{L_{w_{\bal},p }(S)}
&\approx \,\text{inf}\,\big(\Vert  f-g\Vert  _{L_{w_{\bal},p  }(S)}
+\frac 1n\, \Vert  \wt P_{S,\bal}(D)g\Vert  _{L_{w_{\bal},p
}(S)}\big)\\
&\equiv K\big(f,\wt P_{S,\bal }(D), n^{-1}\big)_{L_{w_{\bal},p}(S)}
\end{aligned}
\end{equation}
where $\wt P_{S,\bal} $ is given by (\ref{Eq12.25}).  For
$w_{\bal}(\bx) = 1$ we denote $M_{n,\bal} \equiv M_n$ and $\wt
P_{S,\alpha  }(D) \equiv \wt P_S(D).$  Berens et~al. (see
\cite{Be-Sc-Xu}) conjectured that for $\bxi \in E_S$
\begin{equation}\label{Eq12.28}
\Big\Vert  \pd{}{\bxi  }\;\wt d_S(\bx,\bxi  )\;\pd{}{\bxi  }\; f\Big\Vert
_{L_p(S)} \le C\Vert  \wt P_S(D)f\Vert  _{L_p(S)}\q\text{\rm for} \q 1<p<\infty
\end{equation}
and proved (\ref{Eq12.28}) for $p=2.$  In fact, for $p=2$
(\ref{Eq12.28}) was proved for the weighted case as well (see
\cite{Ch-Di93}).  I believe that (\ref{Eq12.28}) is valid for
$1<p<\infty  $ even for the weighted case.  The inequality
(\ref{Eq12.28}) would have some worthwhile applications if proved.
(For instance, the equivalence between the $K$-functional in
(\ref{Eq12.27}) when $w_{\bal}  =1$ and $\wt \omega
^2_S\big(f,\frac{1}{\sqrt n}\,\big)_p.)$

The sharp Marchaud inequality on the simplex with the $K$-functional
given in (\ref{Eq12.27}) was proved in \cite[Theorem 5.1]{Da-Di07}
and is given by
\begin{equation}\label{Eq12.29}
K_{2\beta  }\big(f,\wt P_{S,\bal}(D)^\beta  ,t^{2\beta  }\big)_p \le
Ct^{2\beta  } \Big\{\int^C_t\;\frac{K_{2\gamma
}\big(f,P_{S,\bal}(D)^\gamma  ,u^{2\gamma  }\big)^q_p}{u^{2\beta
q+1}} \;du\Big\}^{1/q}
\end{equation}
for $1<p<\infty  ,$ $\beta  <\gamma  $ and $q=\min\,(p,2)$ where
$\wt P_{S,\bal}(D)$ is given in (\ref{Eq12.22}) and the
$K$-functionals are defined following (\ref{Eq4.10}).

The Nikol'skii inequality on $I_d = [-1,1]\times \dots \times
[-1,1]$ (see \cite[6.9]{Di-Ti}) is given by
\begin{equation}\label{Eq12.30}
\Vert  P_n\Vert  _{L_{w_{\bal,\bbe},q}(I_d)} \le Cn^{\gamma  (\frac
1p-\frac 1q)}\Vert  P_n\Vert  _{L_{w_{\bal,\bbe},p}(I_d)}
\end{equation}
for $0<p<q\le \infty  ,$ $w_{\bal,\bbe}(\bx) = \os{d}{\us{i=1}\prod}
(1-x_i)^{\alpha  _i}(1+x_i)^{\beta  _i}$  with $\alpha
_i,\beta  _i >-1$
and
\newline $\gamma  =\os
d{\us{i=1}\sum}\,\max\big(2+2\,\max(\alpha  _i,\beta  _i),1\big).$

The results (\ref{Eq12.19}) - (\ref{Eq12.29}) can  easily be extended
to replace $S$ with a cube.  In fact, following remarks in
\cite[Section 5]{Di95I}, these results can be extended to Cartesian
products of simplices.  The multivariate Jackson result for weighted
doubling or for Freud-type weights was not studied.

\section{Ul'yanov-type result}\label{Sec13}
For trigonometric polynomials Ul'yanov established relations between
moduli of smoothness in $L_p(T)$ and moduli of smoothness in
$L_q(T),$ $p<q.$ (For the most general form of these types of
relations on $T^d$ see \cite[Section 2]{Di-Ti05}.)
A Ul'yanov-type inequality shows quantitatively how measures of
smoothness of $f$ in $L_p$ influence the measure of smoothness or
the norm of $f$ in $L_q$ when $q>p.$
Here we present
first the analogous relations for $\omega  ^r_\varphi  (f,t)_p$
proved in \cite{Di-Ti05}.

For $f\in L_p[-1,1],$ $0<p<q\le \infty  ,$ $\omega  ^r_\varphi
(f,t)_p$ and $E_n(f)_p$ given by (\ref{Eq1.1}) and (\ref{Eq2.7})
respectively, we have (see \cite[Section3]{Di-Ti05})
\begin{equation}\label{Eq13.1}
\Vert  f\Vert  _{L_q[-1,1]} \le C\Big[\Big\{\int^1_0 \big(u^{-\theta  }
\omega  ^r_\varphi  (f,u)_p\big)^{q_1}\;\frac{du}{u}\Big\}^{1/q_1} +
\Vert  f\Vert  _{L_p[-1,1]}\Big],
\end{equation}
\begin{equation}\label{Eq13.2}
\omega  ^r_\varphi  (f,t)_q\le C\Big(\int^t_0 \big(u^{-\theta
}\omega  ^r_\varphi  (f,u)_p\big)^{q_1}\;\frac{du}{u}\Big)^{1/q_1},
\end{equation}
\begin{equation}\label{Eq13.3}
\Vert  f\Vert  _{L_q[-1,1]} \le C\Big[\Big\{\sum^\infty  _{k=1}
k^{q_1\theta  -1} E_k(f)^{q_1}_p\Big\}^{1/q_1} + \Vert  f\Vert
_{L_p[-1,1]}\Big],
\end{equation}
and
\begin{equation}\label{Eq13.4}
E_n(f)_q \le C\Big\{\sum^\infty  _{k=n} k^{q_1\theta  -1}
E_k(f)^q_p\Big\}^{1/q_1}
\end{equation}
where
$$
q_1 =\begin{cases} q, &q<\infty  \\ 1, &q=\infty  \end{cases}
\q
\text{\rm and}\q
\theta  = 2\big(\frac 1p - \frac 1q\big).
$$
While (\ref{Eq13.1}) and
(\ref{Eq13.2}) are valid for $r=1,2,\dots$ they are useful only for
$r$ big enough $(r> 2(\frac 1p -\frac 1q)$ for $p\ge 1$ and $r+\frac
1p - 1>2(\frac 1p - \frac 1q)$ for $0<p<1).$

For a simple polytope $S\subset \IR^d$ (see Section \ref{Sec12}) the
inequalities (\ref{Eq13.1})~-~(\ref{Eq13.4}) were generalized in
\cite[Section 8]{Di-Ti05}.  It was proved that for $0<p<q\le \infty  ,$
$\wt \omega  ^r_S(f,t)_{L_p(S)}$ and $E_n(f)_{L_p(S)}$ given by
(\ref{Eq12.6}) and (\ref{Eq12.9}) respectively, one has
\begin{equation}\label{Eq13.5}
\Vert  f\Vert  _{L_q(S)} \le C\Big[\Big\{\int^1_0 \big(u^{-\theta
}\wt \omega
^r_S(f,u)_{L_p(S)}\Big)^{q_1}\;\frac{du}{u}\Big\}^{1/q_1} + \Vert
f\Vert  _{L_p(S)}\Big],
\end{equation}
\begin{equation}\label{Eq13.6}
\wt \omega  \,^r_S(f,t)_{L_q(S)} \le C\Big(\int^t_0 \big(u^{-\theta
} \,\wt\omega
\,^r_S(f,u)_{L_p(S)}\big)^{q_1}\;\frac{du}{u}\Big)^{1/q_1},
\end{equation}
\begin{equation}\label{Eq13.7}
\Vert  f\Vert  _{L_q(S)} \le C\Big[\Big\{\sum^\infty  _{k=1}
k^{q_1\theta  -1} E_k(f)^{q_1}_{L_p(S)}\Big\}^{1/q_1} + \Vert
f\Vert  _{L_p(S)}\Big],
\end{equation}
and
\begin{equation}\label{Eq13.8}
E_n(f)_{L_q(S)} \le C\Big\{\sum^\infty  _{k=n} k^{q_1\theta  -1}
E_k(f)^q_p\Big\}^{1/q_1}
\end{equation}
where
$$
q_1 = \begin{cases} q , &q<\infty  \\ 1, &q=\infty
\end{cases}\q \text{\rm and}\q \theta  = 2d\big(\frac 1p - \frac 1q\big).
$$

In fact, (\ref{Eq13.5}) - (\ref{Eq13.8}) contain (\ref{Eq13.1}) -
(\ref{Eq13.4}) (when $d=1)$ and (\ref{Eq13.1}) - (\ref{Eq13.4}) were
presented here explicitly for those interested mainly in the
one-dimensional case and in the moduli defined by $\omega
^r_\varphi  (f,t)_p,$ which is less intricate than  $\wt\omega
^r_S(f,t)_{L_p(S)}.$

For the weighted $L_p(\IR)$ with the Freud weight $W_\alpha  (x) =
\;\exp(-\vert  x\vert  ^\alpha  )$ a set of Ul'yanov-type inequalities was given
in \cite[Section 9]{Di-Ti05}.  For $0<p<q\le \infty  ,$ $W_\alpha
=\;\exp (-\vert  x\vert  ^\alpha  )$ $(\alpha  >1),$
$E_n(f)_{W_\alpha  ,p}$ and $\omega  ^r(f,W_\alpha  ,t)_p$ given by
(\ref{Eq11.9}) and (\ref{Eq11.5}) respectively, we have
\begin{equation}\label{Eq13.9}
\Vert  W_\alpha  f\Vert  _{L_q(\IR)} \le C\Big[\Big\{\sum^\infty
_{k=1} k^{q_1\theta  -1} E_k(f)^{q_1}_{W_\alpha  ,p}\Big\}^{1/q_1} +
\Vert  W_\alpha  f\Vert  _{L_p(\IR)}\Big],
\end{equation}
\begin{equation}\label{Eq13.10}
E_n(f) _{W_\alpha  ,q} \le C\Big\{\sum^\infty  _{k=n} k^{q_1\theta
-1}E_k(f)^{q_1}_{W_\alpha  ,p}\Big\}^{1/q_1},
\end{equation}
\begin{equation}\label{Eq13.11}
\Vert  W_\alpha  f\Vert  _{L_q(\IR)} \le C\Big[\Big\{\int^1_0
\big(u^{-\eta  }\omega  ^r (f,W_\alpha
,t)_p\big)^{q_1}\;\frac{du}{u}\Big\}^{1/q_1} + \Vert  W_\alpha
f\Vert  _{L_p(\IR)}\Big],
\end{equation}
and
\begin{equation}\label{Eq13.12}
\omega  ^r(f,W_\alpha  ,t)_q \le C\Big\{\int^t_0 \big(u^{-\eta
}\omega  ^r(f,W_\alpha  ,t)_p\big)^{q_1}\;\frac{du}{u}\Big\}^{1/q_1}
\end{equation}
where
$$
q_1 =\begin{cases} q, &q<\infty  \\ 1, &q=\infty
\end{cases},\q \theta   = \frac{\alpha  -1}{\alpha  } \,\big(\frac
1p -\frac 1q\big)\q \text{\rm and} \q \eta  =\frac 1p - \frac 1q\,.
$$
It turns out that the Nikol'ski-type inequalities and realization
results using best approximants following (\ref{Eq5.5}),
(\ref{Eq11.11}) or (\ref{Eq12.18}) are crucial for the proof of the
above-mentioned Ul'yanov-type inequalities.

We note that an inequality like (\ref{Eq13.1}) can be stated using
Besov spaces terminology.  The Besov space $B^\theta
_{p,q}(\varphi,r  )$ is given by the norm or quasi-norm
\begin{equation}\label{Eq13.13}
\Vert  f\Vert  _{B^\theta  _{p,q}(\varphi,r  )} =
\Big(\int^1_0\big(u^{-\theta  }\omega  ^r_\varphi
(f,u)_p\big)^{q}\,\frac{du}{u}\Big)^{1/q} + \Vert  f\Vert
_{L_p[-1,1]}.
\end{equation}

The inequality (\ref{Eq13.1}) means that for $\theta  =2\big(\frac
1p - \frac 1q\big)$ and $0<p<q<\infty  ,$ $B^\theta  _{p,q}(\varphi
,r)$ is continuously embedded in $L_q[-1,1],$ which can be written as
\begin{equation}\label{Eq13.14}
B^\theta  _{p,q}(\varphi  ,r) \hookrightarrow L_q[-1,1].
\end{equation}

In \cite{Di-Ti05} examples are given to show that the power $q_1=q$ is
optimal when $q<\infty  .$

\section{$\omega  ^r_{\varphi  ^\lambda  }(f,t)_\infty  ,$ $0\le
\lambda  \le 1,$ filling the gap}\label{Sec14}
For $C[-1,1]$ and $\omega  ^r(f,t)_{C[-1,1]}$ the classical
``pointwise estimate'' theory established by Dzyadic, Timan, Brudnyi
and others (see \cite{TiA}) yields a complete (pointwise)
description of polynomial approximation on $C[-1,1].$  Estimates
using $\omega  ^r_\varphi  (f,t)_{C[-1,1]}$ yield a complete (norm)
description of polynomial approximation on $C[-1,1].$  While the
estimates using $\omega  ^r_\varphi  (f,t)_p$ are applicable to all
$p,$ $0<p\le \infty  ,$ one has two different ways to characterize
polynomial approximation when $p=\infty  $ (for which the relevant
theory is on $C[-1,1]).$  In an effort to unify these two theories
(for $C[-1,1]),$ one can use the moduli $\omega  _{\varphi ^ \lambda
}(f,t)_{C[-1,1]}$ given for $0\le \lambda  \le 1$ and $\varphi  (x)
= \sqrt{1-x^2}$ by
\begin{equation}\label{Eq14.1}
\omega  ^r_{\varphi  ^\lambda  }(f,t)_{C[-1,1]} = \us{\vert  h\vert
\le t}\sup\;\Vert  \Delta  ^r_{h\varphi  ^\lambda  }f\Vert
_{C[-1,1]}
\end{equation}
where
\begin{equation}\label{Eq14.2}
\Delta  ^r_{h\varphi  ^\lambda  }f(x) = \begin{cases}
\os r{\us{k=0}\sum} (-1)^k\binom rk f\big(x+(\frac r2-k)h\varphi
^\lambda  (x)\big) &\text{\rm when}\q x\pm \frac r2\,h\varphi
^\lambda  (x) \in [-1,1]\\
0 &\text{\rm otherwise}
\end{cases}
\end{equation}
(see \cite{Di-Ji}).

Clearly, $\omega  ^r_{\varphi  ^\lambda  }(f,t)_{C[-1,1]}$ is
$\omega  ^r(f,t)_{C[-1,1]}$ when $\lambda  =0$ and it is $\omega
^r_\varphi  (f,t)_{C[-1,1]}$ when $\lambda  =1.$

The direct estimate using $\omega  ^r_{\varphi  ^\lambda
}(f,t)_{C[-1,1]}$ proved in \cite[Theorem 2.1]{Di-Ji} states that
for
\newline $f\in C[-1,1]$ there exists a sequence of polynomials $P_n$ that
satisfies
\begin{equation}\label{Eq14.3}
\vert  f(x)-P_n(x)\vert  \le C(r,\lambda  )\omega  ^r_{\varphi
^\lambda  }\big(f,n^{-1}\delta  _n(x)^{1-\lambda  }\big)_{C[-1,1]}
\end{equation}
where $\delta  _n(x) =n^{-1} +\sqrt{1-x^2}$ and $C(r,\lambda  )$ is
independent of $f$ and $n.$  The inequality (\ref{Eq14.3}) fills the
gap between (\ref{Eq2.6}) for $C[-1,1]$ (when $\lambda  =1)$ and the
classical estimate (when $\lambda  =0).$  The converse result with
$\omega  ^r_{\varphi  ^\lambda  }(f,t)_{C[-1,1]}$ was given in
\cite[Theorem 5.1]{Di-Ji} as follows.  For $f\in C[-1,1]$ and
$\omega  (t)$ an increasing function satisfying for some $s$
\begin{equation}\label{Eq14.4}
\omega  (\mu  t) \le C(\mu  ^s+1)\omega  (t),
\end{equation}
the existence of a sequence of polynomials $P_n$ satisfying
\begin{equation}\label{Eq14.5}
\vert  f(x) - P_n(x)\vert  \le M\omega  \big(n^{-1}\delta
_n(x)^{1-\lambda  }\big)
\end{equation}
implies
\begin{equation}\label{Eq14.6}
\omega  ^r_{\varphi  ^\lambda  }(f,t) \le Mt^r \sum_{0<n\le 1/t}
n^{r-1}\omega  (n^{-1}).
\end{equation}
Other estimates were given as well, and results using $\omega
^r_{\varphi  ^\lambda  }(f,t)_{C[-1,1]}$ were followed in many
papers which are not referenced here.

The results mentioned in this section, particularly (\ref{Eq14.3})
and (\ref{Eq14.6}), answer a natural question, and filling the gap
was a necessary endeavor.  However, I feel that one is better off
dealing with either $\lambda  =1$ (and the norm estimate) or with
$\lambda  =0$ (and the pointwise estimate).

\section{Shape-preserving polynomial approximation}\label{Sec15}
Sometimes it is desirable that the polynomial approximating a function
on a given interval have the same shape there as the function
itself.  For example, one may want to approximate a nondecreasing or
convex function on $[-1,1]$ by a nondecreasing or convex polynomial
on $[-1,1].$  This aspect of polynomial approximation has attracted
much attention and dozens of papers have been published, mostly in
the last twenty years, covering its many variations.  It is clear to me
that in this survey, I will not be able to do justice to the topic,
which may require a separate survey.  I refer the reader to a survey
by Leviatan (see \cite{Le}) and two subsequent papers by Kopotun,
Leviatan and Shevchuk (see \cite{Ko-Le-Sh05} and \cite{Ko-Le-Sh06})
where many of the related results are described. (The words
``final frontier'' and ``conclusion'' in the last two articles do
not mean that the whole subject of shape-preserving polynomial
approximation
is to be abandoned by these
authors.)  One can probably consider this section as an introduction
to the subject, rather than a survey of the main results.

For a long time it was known that if $f(x)$ satisfies $\Delta  ^k_hf
(x)\ge 0$ on $\big[\frac{kh}{2},\,1-\frac{kh}{2}\big ]$ for some $k$ $(k=0,1,2,\dots),$ then its
Bernstein polynomials, $B_n(f,x)$ given in (\ref{Eq8.1}) satisfy
$\Delta  ^k_h B_n(f,x)\ge 0$ (or $(\frac{d}{dx})^k B_n(f,x)\ge 0)$ for
that $k.$  We recall that $\Delta  ^k_h f(x)> 0$ represents a
condition on the shape of $f;$ for example, when $k=0,$ then $f$ is
positive, when $k=1,$ $f$ is nondecreasing, and when $k=2,$ then $f$ is convex
etc. (Recall $\Delta  ^k_h f(x)$ is given by (\ref{Eq2.3}).)

It is known that
$$
\vert  B_n(f,x) - f(x)\vert  \le C\omega
^2\big(f,\sqrt{\frac{x(1-x)}{n}}\,\big)_{C[0,1]}
$$
(the pointwise
estimate) and
$$\Vert  B_nf-f\Vert  _{C[0,1]} \le C\omega  ^2_\varphi
(f,1/\sqrt n)_{C[0,1]}\q\text{\rm with}\q\varphi  ^2 = x(1-x)
$$
(the norm
estimate).

The approximation by a general polynomial gives rise to
faster convergence, that is $n^{-1/2}$ is replaced by $n^{-1}$ and a
higher degree of smoothness may be considered.  The problem of
shape-preserving polynomial approximation is the relation between
the shape that is preserved and the rate of approximation achievable
under this constraint.

One defines the best constrained polynomial approximation of $f$
satisfying $\Delta  ^k_h f(x)\ge 0$ on $[-1,1]$ by
\begin{equation}\label{Eq15.1}
E^{(k)}_n (f)_{p} = E^{(k)}_n(f)_{L_p[-1,1]} = \;\inf\,\Big(\Vert
f-P_n\Vert  _{L_p[-1,1]}:\Delta  ^k_h P_n(x)\ge 0\;\text{\rm in}\;
[-1,1], \; P_n\in \Pi_n\Big),
\end{equation}
where $\Delta  ^r_hf(x)\ge 0$ in $[-1,1]$ means that for all $x$ and
$h$
\begin{equation}\label{Eq15.2}
\Delta  ^r_h f(x) = \sum^r_{\ell=0} (-1)^\ell \binom r\ell
f\big(x+(\frac r2-\ell)h\big) \ge 0\q\text{\rm where}\q x\pm \frac
{rh}{2} \in [-1,1].
\end{equation}

Shvedov proved (see \cite[Theorem 3]{Sh}) that for any constant $A>0$
and $1\le p\le \infty  $ there exists a function $f\in
C^{(k)}[-1,1]$ such that $\Delta  ^k_hf(x)\ge 0$ in $[-1,1]$ and
\begin{equation}\label{Eq15.3}
E^{(k)}_n (f) _{L_p[-1,1]} \ge A\omega
^{k+2}(f,1/n)_{L_p[-1,1]}\;\text{\rm for}\; n\ge k+2
\end{equation}
where
\begin{equation}\label{Eq15.4}
\omega  ^r(f,h)_{L_p[-1,1]} = \;\sup\, \Vert  \Delta  ^r_h
f(\cdot)\Vert  _{L_p[-1+\frac{rh}{2}\,, 1-\frac{rh}{2}\,]}.
\end{equation}
The inequality (\ref{Eq15.3}) shows that not all Jackson-type
results can be followed.  As $\omega  ^r_\varphi  (f,t)_p\le C\omega
 ^r(f,t)_p,$ Shvedov's negative result applies to $\omega
^r_\varphi  (f,t)_p$ as well, though at the time of publication of
Shvedov's article, estimates of polynomial approximation by $\omega
^r_\varphi  (f,1/n)_p$ were not known.  The knowledge that (as
expected) not all Jackson-type estimates for polynomial
approximation can be followed for shape-preserving polynomial
approximation made the pursuit of the remaining possible estimates
more interesting. Recently, (see \cite[Theorems 1 and 2]{Bo-Pr}) it was
shown (in addition to (\ref{Eq15.3})) that for the function
$$
f(x) =
x^{k-1}_+ =\begin{cases} x^{k-1}, &x\ge 0\\ 0, &x<0\end{cases}
$$
which clearly satisfies $\Delta  ^k_h f(x)\ge 0$ (everywhere) one has
\begin{equation}\label{Eq15.5}
E^{(k)}_n(x^{k-1}_+)_{L_p[-1,1]} \ge \frac{C(k,p)}{n^2} \q\text{\rm
for}\q k>3,
\end{equation}
and, as
\begin{equation}\label{Eq15.6}
\omega  ^3_\varphi  (x^{k-1}_+,t)_p\approx \omega
^3(x^{k-1}_+,t)_p \approx t^{3} \q\text{\rm for} \q k>3,
\end{equation}
it follows that $E^{(k)}_n (x^{k-1}_+)_p\ge C\omega
^3(x^{k-1}_+,1/n)_p$ for $n\ge n_0(C),$ $k>3,$ and $p\le \infty  .$  The
same method (see \cite[Remark 5]{Bo-Pr}) shows that
\begin{equation}\label{Eq15.7}
E^{(3)}_n(x^2_+)_p\ge C\omega  ^3(x^2_+,1/n)_p\;\text{\rm for} \;
n\ge n_0(C)\;\text{\rm and}\; p<\infty  .
\end{equation}

For monotonic functions on $[-1,1]$ satisfying $f\in L_p[-1,1],$ the
Jackson theorem for $0<p\le \infty  $ is given in \cite{De-Le-Yu}
by
\begin{equation}\label{Eq15.8}
E^{(1)}_n (f)_p \le C(p)\omega  ^2_\varphi  (f,1/n)_p\,.
\end{equation}

For convex functions on $[-1,1]$, i.e. when $\Delta  ^2_h f(x)\ge 0,$
it was shown that
\begin{equation}\label{Eq15.10}
E^{(2)}_n(f)_p \le C\omega  ^3_\varphi  (f,1/n)_p, \;\text{\rm
for}\; 0<p\le \infty  .
\end{equation}

In fact, it was known earlier that $E^{(2)}_n(f)_p \le C\omega
^2_\varphi  (f,1/n)_p,$ and it was clear that a gap existed between
that result and (\ref{Eq15.3}).  This gap was closed for $p=\infty
$ by Kopotun (see \cite{Ko94}), and following much of his method, for
$0<p<\infty  $ in \cite{De-Hu-Le}.  Kopotun (see \cite[p.~156]{Ko94}) also
gave the analogue for the pointwise Jackson inequality. That is, he
showed that there exists a sequence of convex polynomials $P_n\in
\Pi_n$ such that
\begin{equation}\label{Eq15.11}
\vert  f(x) - P_n(x)\vert  \le C\omega  ^3\Big(f,\,\frac{1}{n^2}
+\frac 1n\;\sqrt{1-x^2}\,\Big)_{C[-1,1]}\,.
\end{equation}

Recently  Bondarenko (see \cite{Bo}) showed that when $\Delta
^3_h f(x)\ge 0$ in $[-1,1],$ one has
\begin{equation}\label{Eq15.12}
E^{(3)}_n(f)_\infty  \le C\omega  ^3_\varphi  (f,1/n)_\infty  .
\end{equation}
In addition, many other related questions were answered, for
instance, simultaneous approximation of a function and its
derivatives under a shape-preserving constraint or the analogous
pointwise estimate under such constraints. As there were over fifty
articles on the subject of this section, I could not describe all
the results or even just quote them.  (At the beginning of this
section, I already referred to other sources, i.e. \cite{Le},
\cite{Ko-Le-Sh05} and \cite{Ko-Le-Sh06}.) Perhaps I will mention
what might be some  unanswered questions:

\vs\noi
{\bf (I)} Is
\begin{equation}\label{Eq15.13}
E^{(k)}_n (f)_p \le C(p,k)\omega  ^2_\varphi  (f,1/n)_p
\end{equation}
valid for all $k,$ $0<p\le \infty  $ and $n\ge n_0(k,p)?$ (This is
known for $k=1,2,3$ see (\ref{Eq15.8}), (\ref{Eq15.10}) and
(\ref{Eq15.12}).)

\vs\noi
{\bf (II)} Can one obtain the estimate
\begin{equation}\label{Eq15.14}
E^{(3)}_n(f)_\infty  \le C\omega  ^4_\varphi  (f,1/n)_\infty
\;\text{\rm for}\; n\ge n_0?
\end{equation}

\section{Average moduli of smoothness (Ivanov's moduli)}\label{Sec16}
In the text by Sendov and Popov (see \cite{Se-Po}) an alternative to
the moduli of smoothness on $[a,b],$ $T$ or $\IR$ is given and is
called averaged moduli of smoothness. These moduli are defined there
(see \cite[p.~7]{Se-Po}) by
\begin{equation}\label{Eq16.1}
\tau  _k(f,t)_{L_p[a,b]} = \Vert  \omega  _k(f,\cdot;t)\Vert
_{L_p[a,b]}
\end{equation}
for a bounded measurable function $f$ where
\begin{equation}\label{Eq16.2}
\omega  _k(f,x,\delta  ) =\us{\vert  h\vert  \le \delta
}{\sup}\,\Big\{\vert  \vec\Delta  ^k_h f(\zeta  )\vert
:\,\zeta  ,\zeta  +kh\in \Big[x-\frac{k\delta
}{2}\,,x+\frac{k\delta  }{2}\Big]\cap [a,b]\Big\}
\end{equation}
and
\begin{equation}\label{Eq16.3}
\vec\Delta  ^k_h f(x) =\begin{cases}
\os k{\us{\ell=0}\sum}
(-1)^{k-\ell}\binom k\ell f(x+\ell h), &x,x+kh \in [a,b]\\
0, &\text{\rm otherwise}.
\end{cases}
\end{equation}

We note that $\tau  _k(f,t)_{L_p[a,b]}$ given above is not
necessarily finite for all $f\in L_p[a,b].$

In a series of articles (see \cite{Iv}) K.~Ivanov introduced
averaged moduli to deal with algebraic polynomial approximation. The
moduli introduced for $1\le p,$ $q\le \infty  $ are given by (see
\cite[p.~187]{Iv})
\begin{equation}\label{Eq16.4}
\tau  _k\big(f;\psi  (t,\cdot)\big)_{q,p} = \Vert  \omega
_k\big(f;\psi  (t,\cdot)\big)_q\Vert  _p
\end{equation}
where
\begin{equation}\label{Eq16.5}
\omega  _k\big(f,x;\psi  (t,x)\big)_q
= \Big[\frac{1}{2\psi  (t,x)} \;\int^{\psi  (t,x)}_{-\psi  (t,x)}
\,\vert  \vec\Delta  ^k_u f(x)\vert  ^qdu\Big]^{1/q}, \; q<\infty
\end{equation}
and
\begin{equation}
\omega  _k\big(f,x;\psi  (t,x)\big)_\infty   =\sup\,(\vert
\vec\Delta  ^k_h f(x)\vert  ;\vert  h\vert  \le \psi  (t,x)\big).
\end{equation}

The restriction $1\le p,q$ is not necessary, and some results in
case $0<p= q<1$ were discussed in \cite{Ta90}, \cite{Ta91},
\cite{Ta}
and \cite{Di-Hr-Iv}.  We note that
(\ref{Eq16.4}) is finite for any $f\in L_q$ for a fixed $t$ and $p.$
Moreover,
for
$f\in L_p,$ $0<p\le \infty  ,$ $\tau  _k\big(f;\psi
(t,\cdot)\big)_{p,q}$ is finite whenever $q\le p.$

For $\psi  (t,x) = t^2+t\,\sqrt{1-x^2}$ and $[a,b] = [-1,1]$ Ivanov
proved \cite[Corollary 5.2]{Iv} that
\begin{equation}\label{Eq16.7}
\tau  _r\big(f,\psi  (t,\cdot)\big)_{p,p} \approx K_{r,\varphi
}(f,t^r)_p\q\text{\rm for}\q 1\le p\le \infty  .
\end{equation}

One also has
\begin{equation}\label{Eq16.8}
\tau  _r\big(f,\psi  (t,\cdot)\big)_{p,p} \approx \omega  ^r_\varphi
 (f,t)_p\q\text{\rm for}\q 0<p\le \infty  ,
\end{equation}
which, for $0<p<1$ was proved by Tachev (see \cite{Ta}) and also
follows from \cite[Section 7]{Di-Hr-Iv}.
Ivanov also treated the weighted $\tau  $ moduli with weights $w$
and $\psi  $ satisfying some mild conditions (see
\cite[(3.9), (3.10) and (3.11)]{Iv}).

 The moduli $\tau  _r\big(f;\psi  (t,\cdot)\big)_{q,p}$
given by (\ref{Eq16.4}) are a somewhat more cumbersome method to
describe smoothness than $\omega ^r _\varphi  (f,t)_p,$ and their
computation is more difficult. However, they have some advantages.  For
instance, the versatility of having separate $q$ and $p$ may prove
useful. Also in many proofs one resorts to local averages for
obtaining results, and in that direction the averaged moduli may
also be helpful. Many of the results of this paper follow for the
averaged moduli because of (\ref{Eq16.8}), and some were proved
by Ivanov directly for the $\tau  $ moduli independent of
(\ref{Eq16.7}) and (\ref{Eq16.8}).

In conclusion, one should keep in mind the concept given in
(\ref{Eq16.4}) and tools developed by K.~Ivanov for possible use in
polynomial approximation and other problems.

\section{Algebraic addition (Felten's moduli)}\label{Sec17}
In the definition of $\omega  ^r_\varphi  (f,t)_p$ it was clear from
the start that $x\pm h\varphi  (x)$ may not be in the interval
$[-1,1]$ and that $\Delta  ^r_{h\varphi  }f\ne \Delta  _{h\varphi
}(\Delta  ^{r-1}_{h\varphi  }f)$ (where $\Delta  ^r_{h\varphi  }$ is
defined by (\ref{Eq1.2})).  With his goal to alleviate those two
inconveniences (or difficulties) M.~Felten (see [Fe,97,I]
and
[Fe,97,II]) defined the elegant addition
\begin{equation}\label{Eq17.1}
a\oplus b = a\,\sqrt{1-b^2} + b\,\sqrt{1-a^2}\q\text{\rm for}\q
a,b\in [-1,1].
\end{equation}
Felten introduced the difference
\begin{equation}\label{Eq17.2}
_*\Delta  _h f(x) = f(x\oplus h)-f(x), \q _*\Delta  ^r_h f(x) = \;
_*\Delta _h\big(\,_*\Delta  ^{r-1}_h f(x)\big).
\end{equation}
He then dealt with the space $L_{p,\varphi  ^{-1}}[-1,1]$ given by
the norm
\begin{equation}\label{Eq17.3}
\Vert  f\Vert  _{p,\varphi  ^{-1}} = \Big\{\int^1_{-1} \vert
f(x)\vert  ^p\;\frac{dx}{\varphi  (x)}\,\Big\}^{1/p}, \q \Vert
f\Vert  _{\infty  ,\varphi  ^{-1}}\equiv \Vert  f\Vert  _\infty   =
\us{-1\le x\le 1}\sup\;\vert  f(x)\vert
\end{equation}
where $\varphi  (x)^2 = 1-x^2.$ (For $p=\infty  $ he considered
$f\in C[-1,1].)$

He defined the moduli of smoothness
\begin{equation}\label{Eq17.4}
W^r_\varphi  (f,t)_p =\us{0<h\le t}\sup\;\Vert  \,_*\Delta  ^r_h
f\Vert  _{p,\varphi  ^{-1}}.
\end{equation}

Felten proved that $x\oplus h \in [-1,1]$ if $x,h\in [-1,1]$ and
used (\ref{Eq17.2}) for iteration, thus overcoming the inconveniences
mentioned above.

For $W^r_\varphi  (f,t)_p$ Felten proved that
\begin{equation}\label{Eq17.5}
E_n(f)_{p,\varphi  ^{-1}} = O(n^{-\alpha  }) \Longleftrightarrow
W^r_\varphi  (f,t)_p = O(n^{-\alpha  })
\end{equation}
for $0<\alpha  <r$ and $1\le p\le \infty  $ where
\begin{equation}\label{Eq17.6}
E_n(f)_{p,\varphi  ^{-1}} =\us{P\in \Pi_n}\inf\,\Vert  f-P\Vert
_{p,\varphi  ^{-1}}\, .
\end{equation}
Furthermore, he showed that
\begin{equation}\label{Eq17.7}
W^r_\varphi  (f,t)_p\approx \inf\,\big(\Vert  f-g\Vert  _{p,\varphi
^{-1}} + t^r\Vert  D^r g\Vert  _{p,\varphi  ^{-1}}\big)
\end{equation}
where $Dg=\varphi  g^\pr$ and $D^r g = D(D^rg)$ and the infimum is
taken on the class of functions for which $\Vert  D^r g\Vert
_{p,\varphi  ^{-1}}$ is bounded and $g^{(r-1)}$ is locally
absolutely continuous.  (One might as well assume that $g^{(r)}$ is
continuous in $(-1,1)$ with no ill effect on (\ref{Eq17.7}) .)

If not for the fact that the weight in (\ref{Eq17.3}), (\ref{Eq17.6})
and (\ref{Eq17.7}) has to be $(1-x^2)^{-1/2},$ this could have been a
very important development. Unfortunately though, that weight seems
to be crucial (it does not work for the weight 1) and that was
perhaps the justified reason why this direction was not pursued.
Still I feel this was an interesting effort.

\section{Generalized translations}\label{Sec18}
For functions on $T$ the translations $T_tf(x) = f(x+t)$ are
multiplier operators given by
\begin{equation}\label{Eq18.1}
T_t f^\wedge (n) = e^{int}\wh f(n) \q\text{\rm where} \q \wh
g(n) = \frac{1}{2\pi} \, \int^\pi_{-\pi} g(x)e^{-inx}dx.
\end{equation}
As trigonometric polynomial approximation is the model for
investigation of algebraic polynomial approximation in so many
directions, translations using multiplier operators have been
examined for this purpose over the last forty years (see
\cite{Lo-Pe}).  Much work was done by Butzer and mathematicians working with him and under his direction. A survey of those works
including some new results was published in 1992 (see
\cite{Bu-Ja-St}).  The Jacobi translation $\tau  _t$ is given in
\cite{Bu-Ja-St} by
\begin{equation}\label{Eq18.2}
(\tau  _t f)^\wedge (k) = \psi  _k (t) \wh f(k)
\end{equation}
where $\psi_k(t) = \varphi  _k(t)/\varphi  _k(1)$ and
$\varphi  _k(t)$ is given by (\ref{Eq10.5}) and $\wh f(k) =
a_k$ of (\ref{Eq10.6}) for $\alpha  ,\beta  >-1$ as described by
(\ref{Eq10.1}), (\ref{Eq10.5}) and (\ref{Eq10.6}).

The moduli of smoothness were given (see \cite[(3.7),
p.171]{Bu-Ja-St}) by
\begin{equation}\label{Eq18.3}
\omega  ^J_s (f,t)_X = \us{{1-t\le h_j<1}\atop{j=1,\dots,s}}{\sup}
\;\Vert  \Delta  ^J_{h_1}\dots \Delta  ^J_{h_s} f\Vert  _X
\end{equation}
where $X$ is an underlying Banach space (mainly $L_{p,w}[-1,1],$ $1\le p\le
\infty  ,$ with $w=w_{\alpha  ,\beta  }$ of (\ref{Eq10.1})) and
$\Delta  ^J_hf$ is given by
\begin{equation}\label{Eq18.4}
\Delta  ^J_h f(x) = \tau  _h f(x) -f(x).
\end{equation}
Properties of $\omega  ^J_s(f,t)_X$ and their relation to best
weighted algebraic approximation are described in \cite{Bu-Ja-St}.
In that investigation a class of sequences $\{\phi(n)\}$
$(\phi(n)\to 0)$ is given for which $E_n(f),$ $\omega  ^J_s(f,1/n)_X$
and the appropriate $K$-functionals behave like $\phi(n)$ (see for
instance \cite[Theorem 4.1, p.~183]{Bu-Ja-St}).  In fact, the natural
gap between $E_n(f)_X$ and other measures of smoothness is left
bigger than necessary and $E_n(f)_X$ and $\omega  ^J_s(f,t) _X$ are
related via some selected sequences and not by direct and weak
converse inequalities.

The advantage of using generalized translation is that
one has the multiplier operators (\ref{Eq18.2}) which yield
commutativity and other nice properties.  The disadvantages are that
the computation of $\tau  _tf(x)$ and $\omega  ^J_s(f,t)_X$ for a
given $f$ is prohibitive.  In fact, rather than learning about the
behaviour of $E_n(f)_X$ using that of $\omega  ^J_s(f,t)_X,$ it is
actually the behaviour of $\omega  ^J_s(f,t)_X$ that we learn about by using
$E_n(f)_X.$  The rate of convergence of $E_n(f)_X$ now has to be
investigated using other moduli of smoothness.

For $f\in L_p,$
$0<p<1,$ $\tau  _tf$ and $\omega  ^J_s(f,t)_p$ cannot be defined.

M.K. Potapov continued to explore relations between the rate of
approximation by algebraic polynomials and generalized translations,
and has published (together with some coauthors) over twenty
articles on the subject in the last twenty years.  Potapov described
in detail generalized translation as an integral operator for
various situations.  This description is far too long and involved
to give here.  Another feature of Potapov's investigation is that a
relation is given between algebraic polynomial approximation in
$L_{p,w_{\alpha  ,\beta  }}[-1,1]$ and a translation induced by the
weight $w_{\mu  ,\nu  }$ and the differential operator $w_{\mu  ,\nu
 }^{-1} \,\frac{d}{dx}\, (1-x^2) w_{\mu  ,\nu  }\,\frac{d}{dx}\,.$
(Relations are given between the pairs $(\mu  ,\nu  )$ and $(\alpha
,\beta  )$ for which the results are valid.) Potapov and his
coauthors in the situation they investigated achieved a direct
(Jackson-type) and a weak converse result, which is an improvement
on proving that for a sequence $\varphi  (n),$ $\varphi  (n)\to 0,$
that satisfies certain conditions,
$$
E_n(f)_{L_p(w_{\alpha  ,\beta
})}\approx \varphi  (n)\Longleftrightarrow \wt \omega  ^r(f,1/n)_X
\approx \varphi  (n).
$$
Relations with appropriate $K$-functionals
were given but not in the form of the usual equivalence (see
\cite[Theorem 3]{Po01I}).

Clearly, moduli that are defined by
integrals or multipliers cannot be defined for $L_p(w_{\alpha
,\beta  })$ when $0<p<1.$
Also, computation of the behaviour of
these moduli is essentially impossible if one does not use relations
with $E_n(f)_{L_p(w_{\alpha  ,\beta  })}$ and learn about
$E_n(f)_{L_p(w_{\alpha  ,\beta  })}$ by using other moduli.

I refer the reader who is interested in the approach of Potapov and
his coauthors to some of his more recent articles, such as
\cite{Po-Ka}, \cite{Po01I}, \cite{Po01II} and \cite{Po05}.

\section{Lipschitz-type and Besov-type spaces}\label{Sec19}
The Besov-type space that is induced by $\omega  ^r_\varphi
(f,t)_p$ is given by the norm or quasinorm
\begin{equation}\label{Eq19.1}
\Vert  f\Vert  _{B^s_{p,q}(\varphi  ,r)} =
\Big(\int^1_0 u^{-sq}\omega  ^r_\varphi
(f,u)^q_p\,\frac{du}{u}\Big)^{1/q} + \Vert  f\Vert  _{L_p[-1,1]}
\end{equation}
for $0<s,$ $0<p\le \infty  $ and $q<\infty  ,$ and by
\begin{equation}\label{Eq19.2}
\Vert  f\Vert  _{B^s_{p,\infty  }(\varphi  ,r)} =\us
u\sup\;\frac{\omega  ^r_\varphi  (f,u)_p}{u^s} + \Vert  f\Vert
_{L_p[-1,1]}
\end{equation}
for $0<p\le \infty  $ and $q=\infty  .$  The norm or quasi-norm
$\Vert  f\Vert  _{B^s_{p,\infty  }(\varphi  ,r)}$ represents a
Lipschitz-type space.

It was shown in \cite[Corollary 7.2.5]{Di-To87} for $1\le p\le
\infty  $ and in \cite[Theorem 1.1]{Di-Ji-Le} for $0<p<1$ that for
$0<s<r$
\begin{equation}\label{Eq19.3}
\omega  ^r_\varphi  (f,u)_p = O(u^s)\Longleftrightarrow E_n(f)_p =
O(n^{-s}).
\end{equation}
Therefore, for $s<r$
\begin{equation}\label{Eq19.4}
\Vert  f\Vert  _{B^s_{p,\infty }(\varphi  ,r)} \approx \us
n\sup\;n^{-s} E_n(f)_p + \Vert  f\Vert  _{L_p[-1,1]}.
\end{equation}

It was proved for $s< r$ and $1\le p\le\infty  $ (see
\cite[Theorem 2.1]{Di-To88} that
\begin{equation}\label{Eq19.5}
\Vert  f\Vert  _{B^s_{p,q}(\varphi  ,r)} \approx \Big\{\sum^\infty
_{n=1} \, n^{sq-1}E_n(f)^q_p\Big\}^{1/q} + \Vert  f\Vert
_{L_p[-1,1]}.
\end{equation}
This equivalence is valid for $0<p<1$ as well, which now follows
easily the proof in \cite{Di-To88} and the realization results (see
Section \ref{Sec5}).  Similar results are valid for the moduli using Freud
weights $\omega  ^r(f,W,t)_p$ given in (\ref{Eq11.5}) and for moduli
$\wt\omega  ^r_S(f,t)_p$ given by (\ref{Eq12.6}).  For $1\le p\le
\infty  $ one can obtain analogues for the $K$-functionals
$K_{2\alpha  }\big(f,(-\,\frac{d}{dx}\,
(1-x^2)\,\frac{d}{dx})^\alpha  ,t^{2\alpha } \big)_p$ (see
(\ref{Eq4.10})), in which case the result is for $s<2\alpha  .$

We observe that the measure of smoothness that is indicated by
belonging to a given Besov space is not as sharp as the measures
discussed in earlier sections, and that the results described in
this section are mere corollaries of  results described in earlier
sections.
We also note
that while for $1\le p\le \infty  $ one has to avoid
in (\ref{Eq19.1}) $s=r$
(justifiably), for $0<p<1$ we have to avoid $r\le s\le r-1+\frac
1p\,.$

\section{Other methods}\label{Sec20}
To investigate $K$-functionals with step weights such as $\varphi
=\sqrt{1-x^2}$ for example, one can use appropriate transformation
of the variable and study $K$-functionals without step weight or
weights, which is simpler.  Recently, Draganov and Ivanov discussed
this method extensively in a long article (see \cite{Dr-Iv}), and it
seems that this direction of investigation will continue in their
forthcoming papers.  The parts relevant to best algebraic polynomial
approximation are Corollary 5.3 (p.~139) and 8.1 (p.~145) of
\cite{Dr-Iv}.  In \cite[p.~146]{Dr-Iv} transformations related to
Bernstein and Szasz-Mirakian operators on $C[0,1]$ and $C[0,\infty
)$ respectively are discussed.  For $L_p[0,1]$ the Kantorovich and
Durrmeyer operators are discussed in \cite[p.~147]{Dr-Iv}.  Of
course, the $K$-functional of the transformed function is simpler,
but often the difficulty is hidden in the fact that we no longer
deal with the original function but with its transformation.
Nevertheless, there are cases in which this method yields a real
advantage.

In a long paper on the subject the talented M.~Dubiner (see
\cite{Dub}) described smoothness by the rate of local polynomial
approximation.  Using local polynomial approximation to obtain
global polynomial approximation is not new and was used extensively
by many authors.  Dubiner's innovation, however, is that he used the
local polynomial approximation as the basis for his investigation
rather than an intermediate step.

The advantage of the method is that it allegedly yields treatment
for multivariate domains that is not accessible by other methods.
The disadvantages are that for cases which were not investigated
earlier by other methods one cannot estimate the behaviour of such
measures of smoothness.  Another important deficiency is the lack of
converse results.

Moreover, I have difficulty in closing what I
perceive to be many gaps in the proofs and in understanding some of
the concepts.  (While some friends assured me that everything is
okay in the article by Dubiner, those that tried to explain were
stuck like myself.)

Apart from the above and the lack of inverse
theorems, it would be nice if direct and weak converse inequalities
for polynomial approximation on a compact subdomain of $\IR^d$ could
be given.  Without weak converse, this type of result is classical
see (\ref{Eq12.13}).  With weak converse, it is not known for
domains as simple as the unit ball in $\IR^d$ when $d>1.$

Operstein (see \cite{Op}) studied the analogue of the classical
theorems on the rate of pointwise polynomial approximation to
$L_p[-1,1].$  Operstein proved for $f\in L_p[-1,1],$ $\rho_n =
2^{-n}(1-x^2)^{1/2} + 2^{-2n}$ and $\omega  (t)$ satisfying $\omega
(t_1+t_2)\le \omega  (t_1) + \omega  (t_2)$ that there exists a sequence
of polynomials
$P_n~\in~\Pi_{2^n+r-1}$ satisfying
\begin{equation}\label{Eq20.1}
\Big\Vert  \Big\{\,\frac{\Vert  f-P_n\Vert  _{L_p[-1,1]}}{\omega
(\rho_n)} \,\Big\}^\infty  _{n=1}\Big\Vert  _{\ell_p} \le C\Big\Vert
 \,\Big\{\,\frac{\omega  ^r(f,2^{-n})_p}{\omega
(2^{-n})}\,\Big\}^\infty  _{n=1}\Big\Vert  _{\ell_p}
\end{equation}
where
\begin{equation}\label{Eq20.4}
\Vert  \{a_n\}\Vert  _{\ell_p} = \Big\{\sum^\infty  _{n=1} \vert
a_n\vert  ^p\Big\}^{1/p}.
\end{equation}

Furthermore, Operstein showed that if for some sequence $\{P_n\},$
$P_n\in \Pi_{2^n}$ one has
\begin{equation}\label{Eq20.2}
\Big\Vert  \,\Big\{\,\frac{\Vert  f-P_n\Vert  _{L_p[-1,1]}}{\omega
(\rho_n)} \,\Big\}^\infty  _{n=1}\Big\Vert  _{\ell_p} \le 1, \q
1<p<\infty ,
\end{equation}
then
\begin{equation}\label{Eq20.3}
\omega  ^r(f,t)_p \le Ct^r\Big\{\int^1_t\,\Big(\,\frac{\omega
(u)}{u^r}\Big)^q\,\frac{du}{u}\Big\}^{1/q}, \q \frac 1p +\frac 1q =
1.
\end{equation}
One can view (\ref{Eq20.1}) and (\ref{Eq20.3}) as the analogues of
the pointwise direct and converse result.

\section{Epilogue}\label{Sec21}
I have endeavored to mention all directions and progress made
regarding the rate of polynomial approximation in the last twenty
years. Even though my list of references is quite long, there are
over a hundred possible references that might also have been
included.

The issue of best constants was not considered.  Still, I hope that
this survey will be helpful to students and researchers interested
in quantitative estimates of polynomial approximation.

I would like to thank
F. Dai, A.~Prymak and S.~Tikhonov for reading a draft of this
manuscript and eliminating many misprints.

\section{Appendix}

The Jackson-type inequality
$$
E^*_n(f)_B = \inf (\Vert  f-T_n\Vert  _B:\,T_n\in \pmb{\cal T}  _n) \le
C\omega  ^r(f,1/n)_B \eqno{(2.1)^\pr}
$$
for a Banach space $B$ on $T$ satisfying (\ref{Eq2.4}) and
(\ref{Eq2.5}) is essentially known, but as I could not locate a
reference for the exact form $(2.1)^\pr, $ I am adding a proof here.
 One first observes that
$$
\Vert  \sigma  _nf\Vert  _B\le \Vert
f\Vert  _B\q \text{\rm where} \q\sigma  _nf = \frac{1}{2\pi n}\,\int^\pi_{-\pi}\,
\Big(\frac{\sin\,\frac{n(t-x)}{2}}{\sin\,\frac{t-x}{2}}\Big)^2
f(t)dt.
$$

Therefore, $E_{2n}(f)_B \le \Vert  f-2\sigma  _{2n}f +\sigma
_nf\Vert  _B\le 4E_n(f)_B.$  We define $F= f*g =\int^\pi_{-\pi}
f(x+t)g(t)dt$ with $g$ of norm $1$ in $B^*,$ the dual of $B.$  For
appropriately chosen $g$ (and still $\Vert  g\Vert  _{B^*} =1)$
$$
\Vert  f- 2\sigma  _{2n}f+\sigma  _nf\Vert  _B-\varepsilon  \le
\vert  F(0) - 2\sigma  _{2n}F(0) + \sigma  _n F(0)\vert  .
$$
We now have
$$
\begin{aligned}
E_{2n}(f)_B-\varepsilon  &\le \Vert  f-2\sigma  _{2n}f +\sigma  _n
f\Vert  _B-\varepsilon  \\
&\le \vert  F(0) - 2\sigma  _{2n}F(0) + \sigma  _n F(0)\vert  \\
&\le \Vert  F-2\sigma  _{2n} F +\sigma  _n F\Vert  _{C(T)}
\le 4E_n(F)_{C(T)} \\
&\le 4C\omega  ^r(F,1/n)_{C(T)} \le 4C\omega  ^r(f,1/n)_B,
\end{aligned}
$$
and as $\varepsilon  $ is arbitrary (independent of $n)$ and
$2^r \omega  ^r(f,1/2n)_B \ge \omega  ^r(f,1/n)_B,$ $(2.1)^\pr$
is proved.  \qed

\vskip.5in\noi
{Department of Mathematical}
\newline
{and Statistical Sciences}
\newline
{University of Alberta}
\newline
{Edmonton, Alberta}
\newline
{Canada \q T6G 2G1}
\newline
{\tt zditzian@math.ualberta.ca}

\endddoc

\vs
\begin{description}
\item  {2.} Jackson-type estimates
\item  {3.} $K\text{\rm -functionals}$
\item  {4.} $K\text{\rm -functionals} $ (second approach)
\item  {5.} Realization
\item  {6.} Sharp Marchaud and sharp converse inequalities
\item  {7.} Moduli of functions and their derivatives
\item  {8.} Relation with Bernstein polynomial approximation and
other linear operators
\item  {9.} Weighted moduli of smoothness, doubling weights
\item  {10.} Weighted moduli of smoothness, Jacobi weights
\item  {11.} Weighted moduli of smoothness, Freud weights
\item  {12.} Multivariate analogues
\item  {13.} Ul'yanov-type result
\item  {14} $\omega  ^r_{\varphi  ^\lambda  }(f,t)_\infty  ,$ $0\le
\lambda  \le 1,$ filling the gap
\item  {15.} Shape-preserving approximation
\item  {16.} Average moduli of smoothness (Ivanov's moduli)
\item  {17.} Algebraic addition (Felten's moduli)
\item  {18.} Generalized translations
\item  {19.} Other methods
\item  {20.} Conclusions and additional remarks
\end{description}